\documentclass{article}

\newif\ifdraft
\draftfalse 

\usepackage{graphicx} 
\usepackage[T1]{fontenc}
\usepackage{lmodern}
\usepackage{amssymb} 
\usepackage{mathtools} 
\usepackage{amsthm}
\usepackage{bm}
\usepackage{mathrsfs}
\usepackage{etoolbox}
\usepackage{tikz}
\usetikzlibrary{calc,cd}
\usepackage[abbrev]{amsrefs}
\usepackage{comment}

\usepackage{subfiles}

\theoremstyle{plain}
\newtheorem{theorem}{Theorem}[subsection]
\newtheorem{proposition}[theorem]{Proposition}
\newtheorem{lemma}[theorem]{Lemma}
\newtheorem{sublemma}[theorem]{Sublemma}
\newtheorem*{claim*}{Claim}
\newtheorem{corollary}[theorem]{Corollary}
\newtheorem*{theorem*}{Theorem}

\theoremstyle{definition}
\newtheorem{definition}[theorem]{Definition}
\newtheorem{remark}[theorem]{Remark}

\numberwithin{equation}{subsection}

\newcommand{\af}{\mathrm{af}}
\newcommand{\lie}[1]{\mathfrak{#1}}
\newcommand{\affinegcm}[3]{#1_{#2}^{(#3)}}

\newcommand{\Dtwisted}{D_{\ell+1}^{(2)}}
\newcommand{\mixed}{A_{2\ell}^{(2)}}
\newcommand{\bilinearform}[2]{\lparen #1, #2 \rparen}
\newcommand{\pairing}[2]{\langle #1, #2 \rangle}
\newcommand{\realroot}{\Delta_{\mathrm{af}}}
\newcommand{\positiverealroot}{\Delta_{\mathrm{af}}^+}
\newcommand{\cl}[1]{\mathrm{cl}(#1)}
\newcommand{\realcoroot}{\Delta_{\mathrm{af}}^{\vee}}
\newcommand{\positiverealcoroot}{\Delta_{\mathrm{af}}^{\vee,+}}
\newcommand{\poset}[1]{\mathrm{BG}_0(#1)}
\newcommand{\edge}[1]{\xrightarrow{\ #1 \ }}
\newcommand{\LS}[1]{\mathrm{LS}(#1)}
\newcommand{\wt}[1]{\mathrm{wt}(#1)}
\newcommand{\zero}{\boldsymbol{0}}
\newcommand{\realrootsub}[1]{(\Delta_{#1})_{\mathrm{af}}}
\newcommand{\positiverealrootsub}[1]{(\Delta_{#1})_{\mathrm{af}}^+}
\newcommand{\weylsub}[1]{(W_{#1})_{\mathrm{af}}}
\newcommand{\pet}[1]{(W^{#1})_{\mathrm{af}}}
\newcommand{\adjusted}[1]{Q^{#1\text{\textup{-ad}}}}

\newcommand{\sil}{\ell^{\frac{\infty}{2}}}
\newcommand{\sib}{\tfrac{\infty}{2}\mathrm{BG}}

\newcommand{\sils}[1]{\tfrac{\infty}{2}\mathrm{LS}(#1)}

\newcommand{\turn}{\mathrm{Turn(\lambda)}}
\newcommand{\qbg}{\mathrm{QBG}}
\newcommand{\qls}[1]{\mathrm{QLS}(#1)}
\newcommand{\qlsmixed}[1]{\mathrm{QLS}_{\mixed}(#1)}

\DeclarePairedDelimiterX{\set}[1]{\lbrace}{\rbrace}{#1}

\ifdraft
\patchcmd{\section}{\Large}{\large}{}{}
\patchcmd{\subsection}{\large}{\normalsize}{}{}
\fi

\title{Semi-infinite Lakshmibai-Seshadri paths and level-zero extremal weight modules over twisted quantum affine algebras}
\author{Shohei Adachi, Hayato Koike}
\date{}

\begin{document}

\maketitle

\begin{abstract}
    In this paper, we study level-zero extremal weight modules over twisted affine Lie algebras.
    To this end, we introduce semi-infinite Lakshmibai--Seshadri paths associated with a level-zero dominant integral weight $\lambda$.
    We then show that the set $\sils{\lambda}$ of semi-infinite LS paths of shape $\lambda$ is isomorphic, as a crystal, to the crystal basis $\mathcal{B}(\lambda)$ of the corresponding level-zero extremal weight module $V(\lambda)$.
\end{abstract}

\tableofcontents

\section{Introduction}

In \cite{Kashi_modified}, Kashiwara introduced the extremal weight module $V(\lambda)$ over the quantized universal enveloping algebra of a symmetrizable Kac--Moody algebra and also proved the existence of its crystal basis $\mathcal{B}(\lambda)$.
We are interested in the structure of the crystal basis $\mathcal{B}(\lambda)$ of the extremal weight module $V(\lambda)$ when the Lie algebra is of affine type and the extremal weight $\lambda$ is of level-zero.

Let $\mathfrak{g}_\af$ be an affine Kac--Moody Lie algebra over $\mathbb{C}$.
When the affine Lie algebra $\mathfrak{g}_\af$ is of untwisted type, Ishii, Naito, and Sagaki \cite{INS} introduced the semi-infinite Lakshmibai--Seshadri path (in short, LS path) as a combinatorial model for $\mathcal{B}(\lambda)$. In fact, they proved that the set $\sils{\lambda}$ of semi-infinite LS paths of shape $\lambda$ is equipped with a crystal structure and is isomorphic to $\mathcal{B}(\lambda)$ as a crystal for an arbitrary level-zero dominant integral weight $\lambda$ of $\mathfrak{g}_\af$. 
For type $\mixed$, Nomoto\cite{N} introduced the semi-infinite LS path, and showed the same result.

The purpose of this paper is to extend their results to general twisted affine types. 
That is, for the case where $\mathfrak{g}_\af$ is of general twisted type, we introduce a semi-infinite LS path which is a path model equipped with a crystal structure. 
We also prove that $\mathcal{B}(\lambda)$ and the set $\sils{\lambda}$ of semi-infinite LS paths of shape $\lambda$ are isomorphic as crystals for an arbitrary level-zero dominant integral weight $\lambda$ for $\mathfrak{g}_\af$.

Based on the results of this work, it is natural to expect that the graded character of the Demazure submodule can be described by a nonsymmetric Macdonald polynomial.
In this direction, the subset $\mathcal{B}^{\pm}_x(\lambda)$ of $\mathcal{B}(\lambda)$ corresponding to the Demazure submodule $V^{\pm}_x(\lambda)$ should admit a realization as a subset of the set $\sils{\lambda}$ of semi-infinite LS paths.
From this perspective, an interpretation of the Ram--Yip formula in terms of quantum LS paths may offer further insight into the relationship between the graded character and the Macdonald polynomial.
If established, this would yield an analogue of \cite{NNS} in the setting of twisted affine types.

This paper is organized as follows.
In \S2, we fix notation for affine root data.
In \S3, we review the definition of semi-infinite LS paths, the crystal structure on $\sils{\lambda}$, and the proof of the isomorphism theorem for untwisted affine types given in \cite{INS}.
In \S4, building on the material reviewed in \S3, we outline the proof of our result.
We show that it suffices to establish the reduction lemmas (Lemma~\ref{lemma:poset} and Lemma~\ref{lemma:sib}) together with the condition on translation paths (Sublemma~\ref{slem:condition}) in order to prove the isomorphism theorem in our setting.
In \S5, we review the LS paths introduced by Littelmann~\cite{L}, and describe the root data of twisted affine types as well as Peterson's coset representatives.
We also explain the identifications between dual untwisted types and untwisted types, and between type~$\mixed$ and type~$\Dtwisted$.
Moreover, we prove the reduction lemmas, reproduce the propositions in \cite[\S4]{INS}, and equip the set $\sils{\lambda}$ with a crystal structure in the twisted affine cases.
In \S6, we prove Sublemma~\ref{slem:condition} and complete the proof of the isomorphism theorem.
The proof is carried out separately for the cases of dual untwisted types and for type~$\mixed$.
The latter case is more involved, due to the fact that $\pairing{\alpha_{\ell}^{\vee}}{\varpi_{\ell}} = 2$.
The sublemma is proved by restricting the possible forms of positive real roots appearing in $a$-chains satisfying the length condition~\eqref{equation:length_condition}.
This condition arises from the quantum Bruhat graph.
Finally, in \S\ref{appendix:qbg_qls}, we recall the definition of the quantum Bruhat graph and quantum LS paths, and explain their relationship with the semi-infinite Bruhat graph.


\section{Basic notations}

In this section, we review affine root data and also fix our notation.

Let $X_{\af}$ be an affine type generalized Cartan matrix with the index set $I_{\af}$ whose cardinality is $\ell+1$.
Let $\lie{g}_{\af} = \lie{g}(X_{\af})$ be an affine Kac-Moody algebra associated with $X_{\af}$.
Let $\lie{h}_{\af}$ be a Cartan subalgebra of $\lie{g}_{\af}$.
We write the canonical pairing by $\pairing{\cdot}{\cdot}:\lie{h}_{\af} \times \lie{h}_\af^* \rightarrow \mathbb{C}$.
Denote the invariant bilinear form on $\lie{h}_\af^*$ normalized as in \cite[\S2.1]{BN} by $\bilinearform{\cdot}{\cdot}:\lie{h}_\af^* \times \lie{h}_\af^* \rightarrow \mathbb{C}$.
Let $\set{\alpha_i}_{i\in I_\af} \subset \lie{h}_\af^*$ and $ \set{\alpha_i^{\vee}}_{i \in I_{\af}} \subset \lie{h}_\af $ be the set of simple roots and the set of simple coroots of $\lie{g}_{\af}$, respectively. We number the simple roots as in \cite[\S2.1]{BN}.
Then, in the case $X_{\af} = \mixed$, we choose the affine simple root $\alpha_0$ to be the longest one.
We remark that our choice of affine vertex $0$ in type $\mixed$ differs from that in \cite[\S4.8, TABLE Aff 2]{Kac} and \cite{N}.
We set $I \coloneqq I_{\af} \setminus \set{0} = \set{1,\dots,\ell}$.
Let $X$ be the principal submatrix of $X_{\af}$ indexed by $I$; $X$ is a GCM of finite type.
Let $\lie{g} = \lie{g}(X) \subset \lie{g}_{\af}$ be the finite dimensional simple Lie subalgebra of $\lie{g}_{\af}$ associated with $X$, and $\lie{h} \subset \mathfrak{h}_{\af}$ a Cartan subalgebra of $\lie{g}$.
Let $Q \coloneqq\bigoplus_{i \in I}\mathbb{Z}\alpha_i$ and $Q^\vee \coloneqq \bigoplus_{i \in I}\mathbb{Z}\alpha_i^\vee$ be the root lattice and the coroot lattice.
Let $\delta \coloneqq \sum_{i\in I_{\af}}a_i\alpha_i $ and $K \coloneqq \sum_{i \in I_{\af}}a_i^{\vee}\alpha_i^{\vee} $ be the null root and the canonical central element, where $a_i$ and $a_i^{\vee}$ are the labels of $X_{\af}$ and ${}^t\!X_{\af}$, respectively.
Note that $a_0 = 1$ in all cases and
\begin{equation*}
    a_0^{\vee} =
    \begin{cases}
        1 & \text{if $X_{\af} \neq \mixed$,} \\
        2 & \text{if $X_{\af} = \mixed$.}
    \end{cases}
\end{equation*}
We set
\begin{equation*}
    (\mathfrak{h}_{\af}^{*})^{0} \coloneqq \set{\Lambda \in \mathfrak{h}_{\af} \mid \pairing{K}{\Lambda} = 0};
\end{equation*}
we have $(\mathfrak{h}_{\af}^{*})^{0} = \bigoplus_{i \in I}\mathbb{C}\alpha_{i} \oplus \mathbb{C}\delta$.
Let $\mathrm{cl}:(\mathfrak{h}_{\af}^{*})^{0} \twoheadrightarrow (\mathfrak{h}_{\af}^{*})^{0}/\mathbb{C}\delta \cong \mathfrak{h}$ be the canonical projection.
We take an integral weight lattice $P_\af \coloneqq \bigoplus_{i \in I_\af}\mathbb{Z}\Lambda_i \oplus \mathbb{Z}\delta \subset \lie{h}_\af^* $, where $\Lambda_i$ is the $i$-th fundamental weight of $\lie{g}_\af$ for $i \in I_\af$, that is, $\pairing{\alpha_j^{\vee}}{\Lambda_i} = \delta_{ij}$ for all $j \in I_\af$.
For $i \in I$, we define the $i$-th level-zero fundamental weight $\varpi_i \in P_\af$ as follows: when $X_{\af} \neq \mixed$, set
\begin{equation*}
    \varpi_i \coloneqq \Lambda_i - a_i^{\vee}\Lambda_0 \quad \text{for $i \in I$}.
\end{equation*}
When $X_{\af} = \mixed$, set
\begin{gather*}
    \varpi_i \coloneqq
    \begin{cases}
        \Lambda_{i} - \Lambda_{0}     & \text{if $i \neq \ell$,} \\
        2\Lambda_{\ell} - \Lambda_{0} & \text{if $i = \ell$.}
    \end{cases}
\end{gather*}
We remark that $\pairing{K}{\varpi_i} = 0$ for $i \in I$.
We set
\begin{gather*}
    P^0_+ \coloneqq \sum_{i\in I}\mathbb{Z}_{\geq 0}\varpi_i \subset (\mathfrak{h}_{\af}^{*})^{0}.
\end{gather*}
For $\lambda = \sum_{i \in I}m_i\varpi_i \in P^0_+$, we set $J_\lambda \coloneqq \set{i \in I \mid m_i=0} \subset I$.

Let $W_\af \coloneqq \langle s_i \mid i \in I_\af \rangle \subset \mathrm{GL}(\mathfrak{h}_{\af}^{*})$ be the (affine) Weyl group, where $s_i = s_{\alpha_i}$ is the simple reflection with respect to $\alpha_i$. Set $W \coloneqq \langle s_i \mid i\in I\rangle \subset W_\af$, that is, $W$ is the Weyl group of $\lie{g}$.
Let $e \in W_\af$ be the identity element, and $\ell:W_\af \rightarrow \mathbb{Z}_{\ge 0}$ the length function.

Let $\realroot$ be the set of real roots of $\lie{g}_\af$ and $\positiverealroot$ the set of positive real roots of $\lie{g}_\af$.
We set $\Delta \coloneqq \realroot \cap Q$, which is the root system of $\lie{g}$, and $\Delta^+ \coloneqq \Delta \cap \sum_{i\in I}\mathbb{Z}_{\geq 0}\alpha_i$, which is the set of positive roots of $\Delta$.
Denote the reflection with respect to $\beta \in \realroot$ by $s_\beta \in W_\af$.

For $J \subset I$, we set
\begin{gather*}
    Q_J \coloneqq \bigoplus_{j\in J}\mathbb{Z}\alpha_j, \quad Q_J^{\vee} \coloneqq \bigoplus_{j\in J}\mathbb{Z}\alpha_j^{\vee}, \quad Q_J^{\vee,+} \coloneqq \sum_{j\in J}\mathbb{Z}_{\geq 0}\alpha_j^{\vee},\\
    \Delta_J \coloneqq \Delta \cap Q_J, \quad \Delta_J^+ \coloneqq \Delta_J \cap \sum_{i\in I}\mathbb{Z}_{\geq 0}\alpha_i, \quad W_J \coloneqq \langle s_j \mid j \in J \rangle.
\end{gather*}
Let $W^J$ denote the set of minimal-length coset representatives for $W/W_J$.
We know from \cite[\S2.4]{BB} that
\begin{equation*}
    W^J = \set{w \in W \mid w\alpha\in\Delta^+ \text{ for all } \alpha \in \Delta_J^+}.
\end{equation*}
For $w\in W$, $\lfloor w \rfloor = \lfloor w \rfloor^J \in W^J$ denotes the minimal-length coset representative for the coset $wW_J \in W/W_J$.

\section{Review of arguments in \cite{INS}}

In order to explain the differences in construction and proofs in detail, we begin with a brief but sufficiently detailed review of the work of Ishii--Naito--Sagaki \cite{INS}.

\subsection{Definition of semi-infinite LS-paths}
\label{subsection:review_paths}

Let $\mathfrak{g}_\af$ be of untwisted affine type.
For $x\in W_\af$, the \emph{semi-infinite length} of $x=wt_\xi\in W\ltimes Q^\vee$ is defined as
\begin{equation}
    \sil(x) := \ell(w) + 2 \pairing{\xi}{\rho},
\end{equation}
where $\rho \coloneqq (1/2)\sum_{\alpha \in \Delta^{+}}\alpha$.
Let $J \subset I$ and $\pet{J}$ be the set of Peterson's coset representatives (see \cite[\S2.3]{INS} for the definition) for $J$.
We introduce a $\positiverealroot$-labeled directed graph structure on the set $\pet{J}$, where the edges are defined as follows:
\begin{equation*}
    x \edge{\beta} s_{\beta}x
\end{equation*}
such that $s_{\beta}x \in \pet{J}$ and $\sil(s_{\beta}x)= \sil(x) + 1$ with $x \in \pet{J}$ and $\beta \in \positiverealroot$.

We denote this graph by $\sib^{J}$.
\begin{definition}[$a$-chain of shape $\lambda$]
    Let $\lambda \in P^0_+$, $x,y \in \pet{J_{\lambda}}$, and $a \in (0,1]_{\mathbb{Q}} \coloneqq (0,1] \cap \mathbb{Q}$ .
    A directed path
    \begin{equation*}
        x = y_0 \edge{\beta_1} y_1 \edge{\beta_2} \cdots \edge{\beta_{k}} y_k = y
    \end{equation*}
    from $x$ to $y$ in $\sib^{J_\lambda}$ is an $a$-\emph{chain of shape $\lambda$} from $x$ to $y$ if
    \begin{equation*}
        a\pairing{\beta_m^{\vee}}{y_{m-1}\lambda} \in \mathbb{Z}
    \end{equation*}
    for all $1 \leq m \leq k$.
\end{definition}

\begin{definition}[semi-infinite LS paths]
    A \emph{semi-infinite LS path of shape} $\lambda$ is a pair $\eta = (\boldsymbol{x}; \boldsymbol{a}):$
    \begin{equation*}
        \boldsymbol{x}: x_1, \dots , x_s,\quad
        \boldsymbol{a}: 0=a_0 < a_1 < \dots < a_s = 1
    \end{equation*}
    of sequences  in the sets $\pet{J_\lambda}$ and $(0,1]_{\mathbb{Q}}$, respectively, such that
    there is an $a_m$-chain of shape $\lambda$ from $x_{m+1}$ to $x_{m}$ for $1\le m <s$.
    Let $\sils{\lambda}$ denote the set of all semi-infinite LS paths of shape $\lambda$.
\end{definition}

\subsection{Crystal structure on $\sils{\lambda}$}
Fix a level-zero dominant integral weight $\lambda \in P^0_+$.
\begin{theorem}[{{\bf Crystal structure}, \cite[Theorem 3.1.5 (2)]{INS}}]
    \label{thm:A}
    There is a canonical crystal structure on $\sils{\lambda}$.
\end{theorem}
Let $\LS{\lambda}$ be the set of \emph{LS-paths of shape} $\lambda$, which is
defined by Littelmann \cite[\S4]{L} (see \S \ref{subsection:LS_path} for the definition).
In the construction of root operators $e_i,f_i$ ($i \in I_{\af}$) on $\sils{\lambda}$, the following result is crucial.
\begin{lemma}[{{\bf Projection lemma}, \cite[Proposition 3.1.3]{INS}}]
    \label{prop:projection}
    There is a surjective map
    $ \sils{\lambda}\rightarrow{\mathrm{LS}}(\lambda)$
    sending $\eta= (\boldsymbol{x}; \boldsymbol{a})$ to $\overline{\eta}=(x_1\lambda,\ldots,x_s\lambda;a_0,\ldots,a_s).$
\end{lemma}
For the proof of Lemma \ref{prop:projection}, we need some comparison results for graph structures of the level-zero weight poset $\mathrm{BG}_0(\lambda)$ (see \S\ref{subsection:LS_path} for the definition) and $\sib^{J_\lambda}$. For example, we have the following.
\begin{sublemma}[{{\bf Comparison Lemma 1}, \cite[Remark 4.1.3]{INS}}]
    \label{simple_edge}
    Let $\lambda \in P^0_+$.
    For $x \in \pet{J_{\lambda}}$ and $i \in I_\af$,
    \begin{gather*}
        x \edge{\alpha_i} s_ix \text{ in } \sib^{J_\lambda} \Longleftrightarrow \pairing{\alpha_i^\vee}{x\lambda} > 0,\\
        s_ix \edge{\alpha_i} x \text{ in } \sib^{J_\lambda} \Longleftrightarrow \pairing{\alpha_i^\vee}{x\lambda} < 0.
    \end{gather*}
\end{sublemma}
More general comparison result (\cite[Proposition 4.2.1]{INS}) is the following.
\begin{lemma}[{\bf Comparison Lemma 2}]
    \label{Lem:comparison}
    Let $x\in \pet{J_\lambda}$ , and $\beta \in \Delta_{\af}^{+}$. Then, $x\edge{\beta} s_\beta x$ in $\sib^{J_\lambda}$ if and only if $x\lambda \edge{\beta} s_\beta x\lambda$  in $\poset{\lambda}$.
\end{lemma}
This is proved by using ``diamond lemma'' (\cite[Lemma 4.1.4, Lemma 4.1.5]{INS}) for $\sib^{J_{\lambda}}$.
The projection lemma (Lemma \ref{prop:projection}) follows immediately from the comparison lemma 2 (Lemma \ref{Lem:comparison}).

We define the root operators $e_i, f_i$, $i \in I_\af$, on $\sils{\lambda}$ by \emph{lifting} those on $\LS{\lambda}$ (\cite[Definition 3.1.4]{INS}).

In order to finish the proof of Theorem \ref{thm:A}, we need to prove the following.
\begin{lemma}[{{\bf Stability Lemma}, \cite[Theorem 3.1.5 (1)]{INS}}]
    \label{lemma:stability}
    The set $\sils{\lambda}\sqcup\{0\}$ is stable under $e_i,f_i$.
\end{lemma}

This is proved by using the result \cite[Lemma 4.1.6]{INS} concerning $a$-chains in $\sib^{J_{\lambda}}$.

Checking the axioms for crystals is straightforward from the crystal structure on $\LS{\lambda}$, and thus Theorem \ref{thm:A} is proved.
We actually obtain a map of crystals
\begin{equation}
    \sils{\lambda}
    \rightarrow
    \LS{\lambda},
\end{equation}
which is \emph{strict} in the sense of \cite[\S 1.5]{Kashi_modified} (\cite[Remark 4.3.1]{INS}).

\subsection{Isomorphism theorem}

For $\lambda \in P^0_+$, let $V(\lambda)$ denote the extremal weight module of extremal weight $\lambda$ over the quantum affine algebra $U_q(\mathfrak{g}_\af)$(for the detailed definition of extremal weight modules, see \cite[\S3]{Kashi_modified}).
Let $\mathcal{B}(\lambda)$ be a crystal basis of $V(\lambda)$.

\begin{theorem}[{\bf Isomorphism Theorem}]\label{thm:isom}
    There is an isomorphism $\mathcal{B}(\lambda)\cong \sils{\lambda}$ of crystals.
\end{theorem}

Let $u_\lambda \in \mathcal{B}(\lambda)$ be an element corresponding to the generator of $V(\lambda)$, and $\mathcal{B}(\lambda)^{\circ}$ the connected component of $\mathcal{B}(\lambda)$ containing $u_\lambda$.
We use a result due to Beck and Nakajima \cite{BN}, which establishes the isomorphism
\begin{equation}
    \mathcal{B}(\lambda)\cong
    \mathrm{Par}(\lambda)\otimes
    \mathcal{B}(\lambda)^\circ
\end{equation}
of crystals, where $\mathrm{Par}(\lambda)$ is a certain set of
tuples of partitions defined in \cite[(3.2.1)]{INS} (this is ``$\mathbf{N}^{\mathscr{R}_0}(\lambda)'$'' in \cite[\S4]{BN}).
Thus, $\mathrm{Par}(\lambda)$ parametrizes the connected components of $\mathcal{B}(\lambda)$.
We turn to investigate the connected components of $\sils{\lambda}$.
Let $\mathrm{C}(\sils{\lambda})$ denote the set of connected components of $\sils{\lambda}$.
\begin{lemma}[{\bf Connected components of $\sils{\lambda}$, \cite[Proposition 7.2.1]{INS}}]
    \label{lem:conn=conn}
    There exists a bijective correspondence
    \begin{equation*}
        \Theta: \mathrm{Par}(\lambda) \overset{\sim}{\rightarrow} \mathrm{C}(\sils{\lambda})
    \end{equation*}
    between $\mathrm{Par}(\lambda)$ and the set ${\mathrm C}(\sils{\lambda})$  of connected components of $\sils{\lambda}$.
\end{lemma}
In the proof of this lemma, it is essential to prove the following lemma.

\begin{lemma}[{\bf {Special elements}, \cite[Proposition 7.1.2]{INS}}]
    \label{lem:extremal}
    Each connected component of $\sils{\lambda}$ contains a unique element of
    the form
    \begin{equation}
        \label{equation:special}
        (z_{\xi_1}t_{\xi_1},\ldots,z_{\xi_{s-1}}t_{\xi_{s-1}},e;a_0,\ldots,a_s),
    \end{equation}
    with $J_\lambda$-adjusted elements $\xi_1,\dots,\xi_{s-1} \in Q^\vee$ (see \cite[\S2.3]{INS} for the definition).
\end{lemma}

For $\lambda = \sum_{i \in I} m_i\varpi_i \in P^0_+$ and $a \in (0,1]_{\mathbb{Q}}$, we set
\begin{gather}
    J^c(\lambda;a) \coloneqq \set{ i \in J_\lambda^c = I \setminus J_\lambda \mid am_i \in \mathbb{Z}},\\
    I(\lambda; a) \coloneqq \set{ i \in I \mid am_i \in \mathbb{Z} }
    = J_\lambda \sqcup J^c(\lambda; a),
\end{gather}
and set
\begin{gather}
    \turn \coloneqq \left\{ \frac{k}{m_i} \, \middle\vert \, i \in J_\lambda^c, 0 \le k \le m_i \right\}.
\end{gather}

For a subset $J \subset I$, $[\, \cdot \,]_{J^c} : Q^{\vee} = Q_{J^c}^{\vee} \oplus Q_J^{\vee} \twoheadrightarrow Q_{J^c}^{\vee} $ denotes the projection from $Q^{\vee}$ onto $ Q_{J^c}^{\vee} $ with kernel $Q_J^{\vee}$.

Lemma \ref{lem:extremal} is proved by using $N$-multiple map $\sigma_N:\sils{\lambda} \rightarrow \sils{\lambda}^{\otimes N}$ on $\sils{\lambda}$ and the following lemma.

\begin{sublemma}[{\textbf{Condition of translation paths}, \cite[Proposition 7.1.1]{INS}}]
    \label{slem:condition}
    Let $\xi_1, \dots , \xi_s \in Q^{\vee}$ be $J_{\lambda}$-adjusted. An element
    \begin{equation}\label{equation:translation_path_dual}
        \eta = (z_{\xi_1}t_{\xi_1}, \dots, z_{\xi_s}t_{\xi_s}; a_0, \dots, a_s)
    \end{equation}
    is contained in $\sils{\lambda}$ if and only if
    $a_m \in \mathrm{Turn}(\lambda)$ for all $0 \le m \le s$
    and $[\xi_m - \xi_{m+1}]_{J_{\lambda}^{c}} \in Q^{\vee,+}_{J^c(\lambda; a_m)} \setminus \set{0}$ for all $1 \le m \le s-1$.
\end{sublemma}

The proof of Theorem \ref{thm:isom} is completed if the following two results hold.

\begin{lemma}[{{\bf Identification of identity components}, \cite[Proposition 3.2.2]{INS}}]
    \label{lem:id_conns}
    There exists an isomorphism
    \begin{equation}
        \sils{\lambda}^{\circ} \cong \mathcal{B}(\lambda)^\circ
    \end{equation}
    of crystals.
    Here, $\sils{\lambda}^{\circ}$ is the connected component containing $\eta_e \coloneqq (e;0,1) \in \sils{\lambda}$.
\end{lemma}
The proof of this lemma proceeds in two steps.
First, we identify the stabilizer of $u_{\lambda}$ in the action of $W_{\af}$ on $\mathcal{B}(\lambda)$ by using \cite[Remark 4.17]{BN} (see \cite[Proposition 5.1.1]{INS}).
Second, we compare the actions of the root operators on $u_{\lambda} \in \mathcal{B}(\lambda)^{\circ}$ and $\eta_{e} \in \sils{\lambda}^{\circ}$ under two $N$-multiple maps $\sigma_{N}:\mathcal{B}(\lambda)^{\circ} \rightarrow (\mathcal{B}(\lambda)^{\circ})^{\otimes N}$ and $\sigma_N:\sils{\lambda} \rightarrow \sils{\lambda}^{\otimes N}$.

\begin{lemma}[{{\bf Decomposition theorem}, \cite[Proposition 3.2.4]{INS}}]
    \label{lem:fin}
    There exists an isomorphism
    \begin{equation}
        \mathrm{Par}(\lambda)\otimes \sils{\lambda}^\circ \cong \sils{\lambda}
    \end{equation}
    of crystals.
\end{lemma}
Now we have obtained Lemma \ref{lem:conn=conn}, thus what we need to show to prove this lemma is that for each element $\boldsymbol{\rho} \in \mathrm{Par}(\lambda)$, there exists an isomorphism
\begin{equation}
    \label{equation:isom_each_components}
    \Theta(\boldsymbol{\rho}) \cong \set{\boldsymbol{\rho}} \otimes \sils{\lambda}^{\circ}
\end{equation}
of crystals.
The isomorphism is established by comparing how the root operators act on the special element $\eta_{\boldsymbol{\rho}} \in \Theta(\boldsymbol{\rho})$ of the form \eqref{equation:special} and the element $\boldsymbol{\rho} \otimes \eta_{e} \in \set{\boldsymbol{\rho}} \otimes \sils{\lambda}^{\circ}$ under the $N$-multiple map $\sigma_{N}$ on $\sils{\lambda}$.

The proof of Theorem \ref{thm:isom}
goes as follows:
\begin{align*}
    \mathcal{B}(\lambda) &
    \cong
    \mathrm{Par}(\lambda)\otimes
    \mathcal{B}(\lambda)^\circ\quad (\text{Beck-Nakajima})                          \\
                         & \cong
    \mathrm{Par}(\lambda)\otimes
    \sils{\lambda}^\circ\quad\text{(Lemma \ref{lem:id_conns})}                      \\
                         & \cong \sils{\lambda} \quad \text{(Lemma \ref{lem:fin})}.
\end{align*}

\section{Outline of arguments for twisted cases}

To prove the isomorphism theorem (Theorem \ref{thm:A}) in the case of twisted types, we shall follow the logical stream of \cite{INS}.
We discuss the points which we need to consider when we study the twisted cases.

First of all, we need to
define the semi-infinite length on the affine Weyl group $W_\af$
associated with the twisted affine Lie algebra $\mathfrak{g}_\af$.
We know the following identification
\begin{equation*}
    W_\af\cong W\ltimes Q
\end{equation*}
from \cite[Proposition 6.5]{Kac} and the list of the lattices in \cite[\S6.7]{Kac}
(see also \eqref{equation:isom_Weyl_group}).

\begin{definition}[semi-infinite length]
    \label{definition:semi-infinite_length_twisted}
    For $x = wt_\xi \in W_\af \cong W \ltimes Q$ with $w \in W$ and $\xi \in Q$,
    we define the \emph{semi-infinite length} $\sil(x)$ of $x$ by
    \begin{equation}\label{equation:semi-infinite_length}
        \sil(x) := \ell(w) + 2 \pairing{\rho^\vee}{\xi},
    \end{equation}
    where $\rho^\vee := (1/2)\sum_{\alpha \in \Delta^+} \alpha^\vee$.
\end{definition}

Note that when $\mathfrak{g}_{\af}$ is of dual untwisted type, this function coincides with the semi-infinite length on $W_\af$ considered as the Weyl group of the dual Lie algebra $\mathfrak{g}_\af^\lor$ of $\mathfrak{g}_{\af}$.

Once we define the semi-infinite length, the definition of $\sils{\lambda}$ is parallel to the case of untwisted types.

We remark that the construction and property of the $N$-multiple map $\sigma_N:\mathcal{B}(\lambda)^{\circ} \rightarrow (\mathcal{B}(\lambda)^{\circ})^{\otimes N}$ on $\mathcal{B}(\lambda)^{\circ}$ is based on the results by Beck--Nakajima \cite{BN} and Naito--Sagaki \cite{NS_path_model} for general affine types.
Therefore, we can carry over the argument concerning the $N$-multiple map $\sigma_N:\mathcal{B}(\lambda)^{\circ} \rightarrow (\mathcal{B}(\lambda)^{\circ})^{\otimes N}$ in \cite[\S5]{INS} into twisted cases.
In addition, we can also transfer the argument concerning the $N$-multiple map $\sigma_N:\sils{\lambda} \rightarrow \sils{\lambda}^{\otimes N}$ on $\sils{\lambda}$ in \cite[\S5 and \S7]{INS} to our setting because our definition of semi-infinite LS paths is parallel to that in \cite{INS}.
Therefore, in this paper, we omit the arguments regarding $N$-multiple maps.

Based on the review in the previous section and the above note, we present the outline of the proof of our result.

Regarding the crystal structure of $\sils{\lambda}$, if $\mathfrak{g}_\af$ is of dual untwisted type, we reduce the proofs to those for its dual type, and if $\mathfrak{g}_\af$ is of type $A_{2\ell}^{(2)}$, we reduce the proof to that for type $D_{\ell+1}^{(2)}$.
Specifically, we prove the reduction lemma for the level-zero weight poset (Lemma \ref{lemma:poset}) and that for the semi-infinite Bruhat graph (Lemma \ref{lemma:sib}).
From these reduction lemmas, we can reproduce the comparison lemmas and \cite[Lemma 4.1.6]{INS} in the twisted cases.
Therefore, we obtain the projection lemma and the stability lemma.
See \S \ref{subsection:twisted_affine_root_data} and \S \ref{subsection:crystal_structure} for more details.

Regarding the isomorphism theorem, the key results we need to prove are the identification of identity components (Lemma \ref{lem:id_conns}) and the decomposition theorem (Lemma \ref{lem:fin}).
For Lemma \ref{lem:id_conns}, we can prove \cite[Proposition 5.1.1]{INS} in the same way as in the case of untwisted types except that $Q^\vee$ and $Q_J^\vee$  should be replaced with $Q$ and $Q^\vee_J$, respectively.
The remaining part of the proof relies only on the theory of $N$-multiple maps.
Therefore, we can obtain Lemma \ref{lem:id_conns} straightforwardly.
For Lemma \ref{lem:fin}, once we obtain Sublemma \ref{slem:condition}, we can establish the bijection $\Theta:\mathrm{Par}(\lambda) \overset{\sim}{\rightarrow} \mathrm{C}(\sils{\lambda})$ in Lemma \ref{lem:conn=conn} and the isomorphism \eqref{equation:isom_each_components} by using the $N$-multiple map $\sigma_{N}$ on $\sils{\lambda}$.
Hence, we only need to prove Sublemma \ref{slem:condition} for both the cases of dual untwisted types and type $A_{2\ell}^{(2)}$ separately.
Note that in the cases of twisted types, we should replace $Q^{\vee}$ and $Q_{J^{c}(\lambda;a_m)}^{\vee,+}$ in the statement of Sublemma \ref{slem:condition} with $Q$ and $Q_{J^{c}(\lambda;a_m)}^{+}$, respectively.
In particular, the case of type $A_{2\ell}^{(2)}$ is complicated, since $\pairing{\alpha_{\ell}^{\vee}}{\varpi_{\ell}} = 2$.
See \S \ref{section:connected_component} for the proofs.

\section{Semi-infinite Lakshmibai Seshadri paths for twisted affine types}
\label{section:semi-infinite_LS_path}

\subsection{Lakshmibai--Seshadri paths}
\label{subsection:LS_path}

In this subsection, we review the notion of level-zero weight poset and Lakshmibai--Seshadri paths (for short, LS paths).

\begin{definition}[Level-zero weight poset]
    Let $\lambda \in P^0_+$. For $\mu, \nu \in W_\af\lambda$, we write $\mu \le \nu$ if there exist a sequence $\mu = \nu_0, \nu_1,  \dots, \nu_k = \nu$ of elements of $W_\af \lambda$
    and a sequence $\beta_1, \beta_2,  \dots, \beta_k$ of elements of $\Delta^+_\af$
    such that $\nu_m = s_{\beta_m}\nu_{m-1}$ and $\pairing{\beta^{\vee}_m}{\nu_{m-1}} \in \mathbb{Z}_{>0}$ for all $m = 1, 2, \ldots, k$.
    We call the poset $(W_\af \lambda, \le)$ the \emph{level-zero weight poset of shape}   $\lambda$.
    If $\mu \le \nu$, $\operatorname{dist}(\mu,\nu)$ is the maximum length $k$ of all possible such sequences $\mu = \nu_0,\, \nu_1, \, \ldots,\, \nu_k = \nu$.
\end{definition}

\begin{remark}
    The relation $\le$ on $W_{\af}\lambda$ is a partial order on $W_{\af}\lambda$.
\end{remark}

\begin{definition}
    Let $\lambda \in P^0_+$.
    \begin{enumerate}
        \item Let $\poset{\lambda}$ be the graph with vertex set $W_\af\lambda$ and directed edges labeled by positive real roots of the following form:
              \begin{equation*}
                  \mu \edge{\beta} \nu
              \end{equation*}
              such that  $\nu = s_\beta\mu$, $\pairing{\beta^\vee}{\mu} > 0$, and $\operatorname{dist}(\mu,\nu)=1$.
        \item Let $a \in (0,1]_{\mathbb{Q}}$.
              A directed path
              \begin{equation*}
                  \mu = \nu_0 \edge{\beta_1} \nu_1 \edge{\beta_2} \cdots \edge{\beta_k} \nu_k = \nu
              \end{equation*}
              in $\poset{\lambda}$ is called an $a$-\emph{chain} from $\mu$ to $\nu$ if $a\pairing{\beta_m^{\vee}}{\nu_{m-1}} \in \mathbb{Z}$ for all $1 \leq m \leq k$.
    \end{enumerate}
\end{definition}

\begin{definition}[LS path]
    Let $\lambda \in P^0_+$.
    A pair of sequences $\pi = (\nu_1 > \cdots > \nu_s;\, 0 = a_0 < \dots < a_s = 1)$, where $\nu_1, \dots ,\nu_s \in W_\af\lambda$ and $a_0, a_1, \dots ,a_s \in \mathbb{Q}$, is an \emph{LS path of shape} $\lambda$ if there exists an $a_m$-chain from $\nu_{m+1}$ to $\nu_m$ in $\mathrm{BG}_0(\lambda)$ for all $1 \leq m \leq s-1$. Let $\LS{\lambda}$ be the set of LS paths of shape $\lambda$.
\end{definition}

We describe the crystal structure on $\LS{\lambda}$ for $\lambda \in P^0_+$.
We define the map $\mathrm{wt}:\LS{\lambda} \rightarrow P_\af$ by $\wt{\pi} := \sum_{m=1}^{s}(a_{m} - a_{m-1})\nu_{m} \in P_\af$ for $\pi \in \LS{\lambda}$.

For $\pi \in \LS{\lambda}$ and $i \in I_\af$, define the function $H^\pi_i(t)$ on $[0,1]$ and its minimal value  by
\begin{gather*}
    H^\pi_i(t) := \pairing{\alpha^\vee_i}{\pi(t)} \text{ for $t \in [0,1]$,} \\
    m^\pi_i := \min\{ H^\pi_i(t) \mid t \in [0,1] \}.
\end{gather*}
We define the root operators $e_i, f_i$ on the set $\LS{\lambda}$ as follows.
\begin{definition}[root operator]
    Let $\pi = ( \nu_1, \dots, \nu_s ; a_0, \dots, a_s ) \in \LS{\lambda}$.
    Let $\zero$ denote the formal element not contained in $\LS{\lambda}$.
    \begin{enumerate}
        \item If $m^{\pi}_i = 0$, then $e_i \pi := \zero$.
              If $m^{\pi}_i < 0$, then we set
              \begin{equation}\label{def_crystal_for_LS_e}
                  \begin{split}
                      t_1 := & \min \set{ t \in [0,1] \mid H^\pi_i(t) = m^\pi_i },       \\
                      t_0 := & \max \set{ t \in [0,t_1] \mid H^\pi_i(t) = m^\pi_i + 1 }.
                  \end{split}
              \end{equation}
              It is known that $H^\pi_i(t)$ is strictly decreasing on $[t_0, t_1]$ from \cite[Remark~2.2.5]{INS} and \cite[Lemma~4.5 d)]{L}.
              Also, by the minimality of $t_1$, there exists $1 \le q \le s$ such that $t_1 = a_q$.
              Take $1 \le p \le q$ such that $a_{p-1} \leq t_0 < a_p$.
              Then, we define $e_i \pi$ by
              \begin{equation*}
                  \begin{split}
                      e_i \pi = ( & \nu_1, \dots, \nu_{p}, s_i\nu_{p}, \dots, s_i\nu_q, \nu_{q+1}, \dots, \nu_s ; \\
                                  & a_0, \dots, a_{p-1}, t_0, a_p, \dots, a_q = t_1, \dots, a_s).
                  \end{split}
              \end{equation*}
              If $a_{p-1} = t_0$, we drop $\nu_p$ and $a_{p-1}$;
              if $s_i\nu_q = \nu_{q+1}$, we drop $\nu_{q+1}$ and $a_q = t_1$.
        \item If $m^\pi_i - H^{\pi}_i(1) = 0$, then $f_i\pi := \zero$.
              Otherwise, we set
              \begin{equation} \label{def_crystal_for_LS_f}
                  \begin{split}
                       & t_0 := \max \set{ t \in [0,1] \mid H^\pi_i(t) = m^\pi_i },       \\
                       & t_1 := \min \set{ t \in [t_0,1] \mid H^\pi_i(t) = m^\pi_i + 1 }.
                  \end{split}
              \end{equation}
              It is known that $H^\pi_i(t)$ is strictly increasing on $[t_0, t_1]$ from \cite[Remark2.2.5]{INS},\cite[Lemma~4.5 d)]{L}).
              By the maximality of $t_0$, there exists $0 \le p \le s-1$ such that $t_0 = a_p$.
              Take $p \le q \le s-1$ such that $a_{q} < t_1 \le a_{q+1}$.
              Then, we define $f_i\pi$ by
              \begin{equation*}
                  \begin{split}
                      f_i \pi = ( & \nu_1, \dots, \nu_{p},
                      s_i\nu_{p+1}, \dots, s_i\nu_{q+1}, \nu_{q+1}, \dots, \nu_s ;               \\
                                  & a_0, \dots, a_p=t_0, \dots, a_{q}, t_1, a_{q+1} \dots, a_s),
                  \end{split}
              \end{equation*}
              If $a_{q+1} = t_1$, we drop $\nu_{q+1}$ and $a_{q+1}$;
              if $s_i\nu_{p+1} = \nu_q$, we drop $\nu_p$ and $a_p=t_0$.
        \item We set $e_i\zero = f_i\zero := \zero$ for all $i \in I_\af$.
    \end{enumerate}
\end{definition}

For $\pi \in \LS{\lambda}$ and $i\in I_\af$, we set
\begin{gather*}
    \varepsilon_i(\pi) := \max\set{n \in \mathbb{Z}_{\geq 0} \mid e_i^n\pi \neq \zero },\\
    \varphi_i(\pi) := \max\set{n \in \mathbb{Z}_{\geq 0} \mid f_i^n\pi \neq \zero }.
\end{gather*}

\begin{proposition}[see {\cite[\S2 and \S4]{L}}]
    $\LS{\lambda}$ is a crystal whose weight lattice $P_\af$ equipped with the maps $\mathrm{wt}$, $e_i$, $f_i$, $\varepsilon_i$, $\varphi_i$ $(i\in I_\af)$.
\end{proposition}

\subsection{Twisted affine root data and Peterson's coset representatives}
\label{subsection:twisted_affine_root_data}

Suppose that $\mathfrak{g}_{\af}$ is of twisted type.
We define a positive root $\theta_{s} \in \Delta^{+}$ by
\begin{equation*}
    \theta_{s} \coloneqq
    \begin{cases}
        \delta -\alpha_{0}             & \text{if $\mathfrak{g}_{\af}$ is of dual untwisted type}, \\
        \dfrac{\delta - \alpha_{0}}{2} & \text{if $\mathfrak{g}_{\af}$ is of type $\mixed$.}
    \end{cases}
\end{equation*}
This is the highest short root of $\Delta$.
For $\xi \in Q$, denote by $t_\xi \in W_\af$ the translation element with respect to $\xi$ (see \cite[\S6.5]{Kac} for the definition).
We have an isomorphism of groups
\begin{equation}
    \label{equation:isom_Weyl_group}
    \begin{split}
        W_{\af} & \overset{\sim}{\rightarrow} W \ltimes Q \\
        s_{i}   & \mapsto s_{i} \quad \text{($i \in I$)}, \\
        s_{0}   & \mapsto s_{\theta_{s}}t_{-\theta_{s}}.
    \end{split}
\end{equation}

For $\alpha \in \Delta$, we set
\begin{equation*}
    c_\alpha \coloneqq \max\biggl\{1, \frac{\bilinearform{\alpha}{\alpha}}{2} \biggr\},
\end{equation*}
and $c_i \coloneqq c_{\alpha_i}$ for $i \in I$.
Note that $c_{w\alpha} = c_{\alpha}$ for $w \in W$ and $\alpha \in \Delta$.
For $\beta \in \realroot$, $\beta^{\vee}$ denotes the coroot of $\beta$.
We set
\begin{gather*}
    \realcoroot \coloneqq \set{\beta^{\vee} \mid \beta \in \realroot},\\
    \positiverealcoroot \coloneqq \set{\beta^{\vee} \mid \beta \in \positiverealroot}.
\end{gather*}

\subsubsection{Dual untwisted types.}
\label{subsubsection:dual_untwisted}

Suppose $\mathfrak{g}_{\af}$ is of \emph{dual untwisted type}.
Namely, $X_{\af} = \affinegcm{D}{\ell+1}{2}, \affinegcm{A}{2\ell-1}{2}, \affinegcm{E}{6}{2}, \affinegcm{D}{4}{3}$.
For $\alpha \in \Delta$, we have
\begin{equation*}
    c_\alpha =
    \begin{cases}
        1 & \text{if $\alpha$ is a short root in $\Delta$,}                                       \\
        2 & \text{if $X_{\af} \neq \affinegcm{D}{4}{3}$ and $\alpha$ is a long root in $\Delta$,} \\
        3 & \text{if $X_{\af} = \affinegcm{D}{4}{3}$ and $\alpha$ is a long root in $\Delta$.}
    \end{cases}
\end{equation*}
We know that
\begin{gather*}
    \Delta_{\af} = \set{\alpha+nc_\alpha\delta \mid \alpha \in \Delta, n \in \mathbb{Z}},\\
    \positiverealroot = \Delta^+ \sqcup \set{\alpha + nc_\alpha\delta \mid \alpha \in \Delta, n \in \mathbb{Z}_{>0}}.
\end{gather*}

We have
\begin{equation}
    \label{equation:semiproduct_dual}
    s_{\alpha + nc_\alpha\delta} = s_\alpha t_{n\alpha}
\end{equation}
for $\alpha \in \Delta$ and $n \in \mathbb{Z}$.

In the case of dual untwisted types, we consider the affine Lie algebra $\lie{g}_\af^\vee$ whose root data are dual to those of $\lie{g}_\af$.
That is, $\lie{g}_\af^\vee = \lie{g}({}^t\!X_{\af})$ is an untwisted affine Lie algebra.
We identify the objects related to $\lie{g}_\af^\vee$ with the dual objects related to $\lie{g}_\af$.
For example, $\lie{h}_\af^*$ is a Cartan subalgebra of $\lie{g}_\af^\vee$, $\set{\alpha_i^{\vee}}_{i\in I_\af}$ is the set of simple roots of $\lie{g}_\af^{\vee}$, $\set{\alpha_i}_{i \in I_\af}$ is the set of simple coroots of $\lie{g}_\af^{\vee}$, $K$ is the null root of $\lie{g}_\af^\vee$, $\realcoroot$ is the set of real roots of $\lie{g}_\af^\vee$.
Note that the highest short root $\theta_{s}$ of $\mathfrak{g}_{\af}$ is regarded as the coroot of the highest root of $\mathfrak{g}_{\af}^{\vee}$ through this identification.
$W_\af^{\vee} \coloneqq \langle s_{i} = s_{\alpha_{i}^{\vee}} \mid i \in I_{\af} \rangle \subset \mathrm{GL}(\mathfrak{h}_{\af})$ denotes the (affine) Weyl group of $\lie{g}_\af^{\vee}$.
There is an isomorphism of groups $W_{\af}^{\vee} \overset{\sim}{\rightarrow} W \ltimes Q$, given by $s_{i} \mapsto s_{i}$ for $i \in I$ and $s_{0} \mapsto s_{\theta_{s}}t_{-\theta_{s}}$.
Thus, we obtain an isomorphism of groups between $W_\af$ and $W_\af^\vee$ such that
\begin{equation}
    \begin{split}
        \phi: W_\af & \rightarrow W_\af^\vee           \\
        s_i         & \mapsto s_i \quad (i \in I_\af).
    \end{split}
\end{equation}
We remark that
\begin{equation}\label{equation:isom_for_dual_untwisted}
    \phi(x)\beta^\vee = (x\beta)^\vee
\end{equation}
for $x \in W_\af$ and $\beta \in \realroot$.

Let $P_\af^{\vee} \coloneqq \bigoplus_{i \in I_\af} \mathbb{Z}\Lambda_i^\vee \oplus \mathbb{Z}K \subset \lie{h}_\af$ be an integral coweight lattice of $\lie{g}_\af$ identified with an integral weight lattice of $\lie{g}_\af^{\vee}$.
Here, $\Lambda_i^{\vee} \in P_\af^{\vee}$ is the $i$-th fundamental coweight of $\lie{g}_\af$, that is, $\pairing{\Lambda_i^{\vee}}{\alpha_j} = \delta_{ij}$ for all $j \in I_\af$.
For each $i \in I$, we define the $i$-th level-zero fundamental coweight $\varpi_i^{\vee} \coloneqq \Lambda_i^{\vee} - a_i\Lambda_0^{\vee}$.
Note that $\pairing{\varpi_i^{\vee}}{\delta} = 0$ for all $i \in I$.
Set
\begin{equation*}
    P^{0,\vee}_+ \coloneqq \sum_{i \in I} \mathbb{Z}_{\geq 0}\varpi_i^{\vee}.
\end{equation*}

\subsubsection{Type $\mixed$.}

Suppose that $\mathfrak{g}_{\af}$ is of type $\mixed$.
We know that
\begin{gather*}
    \begin{split}
        \Delta_{\af} & = \set{\alpha+n\delta \mid \alpha \in \Delta, n \in \mathbb{Z}}                                                          \\
                     & \hphantom{{}={}}\sqcup\set{2\alpha + (2n-1)\delta \mid \text{$\alpha$ is a short root in $\Delta$, $n \in \mathbb{Z}$}},
    \end{split}\\
    \begin{split}
        \positiverealroot & = \Delta^+ \sqcup \set{\alpha + n\delta \mid \alpha \in \Delta, n \in \mathbb{Z}_{>0}}                                        \\
                          & \hphantom{{}={}}\sqcup\set{2\alpha + (2n-1)\delta \mid \text{$\alpha$ is a short root in $\Delta$, $n \in \mathbb{Z}_{>0}$}}.
    \end{split}
\end{gather*}

We have
\begin{gather}
    s_{\alpha + n\delta} =
    \begin{cases}
        s_\alpha t_{n\alpha}  & \text{if $\alpha$ is a long root in $\Delta$},  \\
        s_\alpha t_{2n\alpha} & \text{if $\alpha$ is a short root in $\Delta$},
    \end{cases}
    \label{equation:semiproduct_mixed1}\\
    s_{2\alpha+(2n-1)\delta} = s_\alpha t_{(2n-1)\alpha}
    \label{equation:semiproduct_mixed2}
\end{gather}
for $\alpha \in \Delta$ and $n \in \mathbb{Z}$.

In the case of type $\mixed$, we consider the affine Lie algebra $\tilde{\lie{g}}_\af = \lie{g}(\Dtwisted)$ of type $\Dtwisted$.
In this case, we follow the notation in \S2, adding a tilde to symbols associated with $\tilde{\mathfrak{g}}_{\af} = \mathfrak{g}(\Dtwisted)$.
For example, $\tilde{\varpi}_{i}$ denotes the $i$-th level-zero fundamental weight of $\tilde{\lie{g}}_{\af}$, and $\tilde{P}^{0}_{+}$ denotes the set of level-zero dominant weights of $\tilde{\mathfrak{g}}_{\af}$.
Note that type $\mixed$ and type $\Dtwisted$ have the same Coxeter diagram and the isomorphic affine Weyl groups.
We define a map $\psi:\Delta_{\af} \rightarrow \tilde{\Delta}_\af$ by
\begin{gather*}
    \psi(\alpha + n\delta) = \alpha + 2n\tilde{\delta}, \\
    \psi(2\alpha + (2n-1)\delta) = \alpha + (2n-1)\tilde{\delta}.
\end{gather*}
Here, we identify $\Delta \subset Q$ and $\tilde{\Delta} \subset \tilde{Q}$.
Note that via this identification of the two root lattices, we have
\begin{equation*}
    2\bilinearform{\xi}{\zeta} = \widetilde{\bilinearform{\xi}{\zeta}}
\end{equation*}
for $\xi, \zeta \in Q$.
Here, $\widetilde{(\cdot,\cdot)}$ denotes the invariant bilinear form on $\tilde{\mathfrak{h}}^{*}_{\af}$, which is the dual space of the Cartan subalgebra of $\tilde{\mathfrak{g}}_{\af} = \mathfrak{g}(\Dtwisted)$.
There exists an group isomorphism given by
\begin{equation}\label{equation:isom_mixed}
    \begin{split}
        \psi:W_\af & \rightarrow \tilde{W}_\af                  \\
        s_i        & \mapsto s_i  \quad \text{($i \in I_\af$).}
    \end{split}
\end{equation}
We remark that
\begin{equation}\label{equation:isom_for_mixed}
    \psi(x)\psi(\beta) = \psi(x\beta)
\end{equation}
for $x \in W_\af$ and $\beta \in \realroot$.

\begin{remark}
    The idea of the identification between type $\mixed$ and type $\Dtwisted$ is based on \cite[\S2]{N}.
    However, note that our choice of the affine node $0$ in the Dynkin diagram of type $A_{2\ell}^{(2)}$ is opposite to that in \cite{N}.
\end{remark}

\subsubsection{Peterson's coset representatives}

Supoose that $\mathfrak{g}_{\af}$ is of general twisted affine type.
Let us describe Peterson's coset representatives.
For $J \subset I$, we set
\begin{gather}
    \realrootsub{J} \coloneqq \set{ \beta \in \realroot \mid \cl{\beta} \in Q_J},\\
    \positiverealrootsub{J} \coloneqq \realrootsub{J} \cap \positiverealroot,\\
    \weylsub{J} \coloneqq W_J \ltimes \set{t_\xi \mid \xi \in Q_J^{\vee}} = \langle s_\beta \mid \beta \in \positiverealrootsub{J} \rangle,\\
    \pet{J} \coloneqq \set{x \in W_\af \mid \text{$x\beta \in \positiverealroot$ for all $\beta \in \positiverealrootsub{J}$}}.
    \label{equation:dual_untwisted_representatives}
\end{gather}

\begin{lemma}\label{lemma:coset}
    \begin{enumerate}
        \item Let $\mathfrak{g}_{\af}$ be of dual untwisted type.
              The group isomorphism $\phi:W_\af \rightarrow W_\af^{\vee}$ maps $\pet{J}$ to $\pet{J}^\vee$, where $\pet{J}^\vee$ is the set of Peterson's coset representatives for $W_\af^\vee$ and $J \subset I$.
        \item Let $\mathfrak{g}_{\af}$ be of type $\mixed$.
              The group isomorphism $\psi:W_\af \rightarrow \tilde{W}_{\af}$ maps $\pet{J}$ to $(\tilde{W}^{J})_{\af}$, where $(\tilde{W}^{J})_{\af}$ is the set of Peterson's coset representatives for $\tilde{W}_{\af}$ and $J \subset I$.
    \end{enumerate}
\end{lemma}
\begin{proof}
    It follows from \eqref{equation:isom_for_dual_untwisted} and \eqref{equation:isom_for_mixed}.
\end{proof}

From this lemma and \cite[Lecture 13]{P} (see also \cite[Lemma 10.6]{LS}), we obtain the following proposition.
\begin{proposition}\label{proposition:unique_factorization_dual_untwisted}
    For all $x \in W_\af$, there exist a unique element $x_1 \in \pet{J}$ and a unique element $x_2 \in (W_J)_\af$ such that $x = x_1x_2$.
\end{proposition}
Define a map $\Pi^J:W_\af \rightarrow \pet{J}$ by $\Pi^J(x) \coloneqq x_1$ in the notation of Proposition~\ref{proposition:unique_factorization_dual_untwisted}.

\begin{definition}
    We say that an element $\xi \in Q$ is \textit{$J$-adjusted} if $\pairing{\alpha^\vee}{\xi} \in \set{-1,0}$ for all $\alpha \in \Delta_J^+$.
    We denote by $\adjusted{J}$ the set of $J$-adjusted elements in $Q$.
\end{definition}

Through the identification explained in the previous subsubsections, we obtain the following lemma from \cite[Lemma 2.3.5]{INS}.

\begin{lemma}\label{lemma:adjusted}
    \begin{enumerate}
        \item For all $\xi \in Q$, there exists a unique element $\phi_J(\xi) \in Q_J$ such that $\xi + \phi_J(\xi) \in \adjusted{J}$. In particular, $\xi \in Q$ is $J$-adjusted if and only if $\phi_J(\xi) = 0$.
        \item For all $\xi \in Q$, there exists an element $z_\xi \in W_J$ such that $\Pi^J(t_\xi) = z_\xi t_{\xi +\phi_J(\xi)}$.
        \item We have
              \begin{equation}
                  \pet{J} = \set{wz_\xi t_\xi \mid \text{ $w \in W^J,\, \xi \in \adjusted{J}$}}.
              \end{equation}
    \end{enumerate}
\end{lemma}

For a subset $J \subset I$, $[\, \cdot \,]_{J^c} : Q = Q_{J^c} \oplus Q_J^{\vee} \twoheadrightarrow Q_{J^c} $ denotes the projection from $Q$ onto $ Q_{J^c} $ with kernel $Q_J$.

\begin{corollary}
    \label{corollary:adjusted}
    Let $J \subset I$ and $\xi,\zeta \in \adjusted{J}$.
    If $[\xi]_{J^{c}} = [\zeta]_{J^{c}}$, then $\xi = \zeta$
\end{corollary}
\begin{proof}
    Write $\xi$ as $\xi = \zeta + (\xi-\zeta)$.
    From $\xi \in \adjusted{J}$ and Lemma \ref{lemma:adjusted}(1), we see that $\phi_{J}(\zeta) = \xi - \zeta$.
    From $\zeta \in \adjusted{J}$ and Lemma \ref{lemma:adjusted}(1), we have $\xi -\zeta = \phi_{J}(\zeta) = 0$.
\end{proof}

\subsection{Semi-infinite Lakshmibai-Seshadri paths}
\label{subsection:semi-infnite_LS_path}

Suppose $\mathfrak{g}_{\af}$ is of twisted affine type.
Recall that we defined the semi-infinite length function for twisted types in Definition \ref{definition:semi-infinite_length_twisted}.

\begin{remark}\label{remark:semi-length_dual}
    When $\mathfrak{g}_\af$ is of dual untwisted type, we have $\sil(x) = \sil(\phi(x))$ for $x \in W_\af$.

    When $\mathfrak{g}_\af$ is of type $\mixed$, we have $\sil(x) = \sil(\psi(x))$ for $x \in W_\af$.
\end{remark}

\begin{remark}
    For $\mixed$ type, Nomoto has already defined the same semi-infinite length function (see \cite[Definition~A.2.1]{N}).
\end{remark}

Now, we can define the semi-infinite Bruhat graph, $a$-chains of shape $\lambda$, semi-infinite LS paths for twisted types in the same way as in \S\ref{subsection:review_paths}.

\subsection{Crystal structure of $\sils{\lambda}$}
\label{subsection:crystal_structure}

In the remainder of this section, we explain the proof of Theorem \ref{thm:A} for twisted types.
As mentioned in the previous section, the proof for the dual untwisted types is reduced to that for the untwisted types, and the proof for type $\mixed$ is reduced to that for type $\Dtwisted$, respectively.
The following lemma is the key to the reduction.

\begin{lemma}\label{lemma:poset}
    \begin{enumerate}
        \item Let $\mathfrak{g}_{\af}$ be of dual untwisted type.
              Let $\lambda = \sum_{i \in I}m_i\varpi_i \in  P^0_+$ and $\lambda^\vee = \sum_{i \in I}m_i\varpi_i^\vee \in P^{0,\vee}_{+}$.
              Let $x \in W_\af$ and $\beta \in \positiverealroot$.
              $x\lambda \edge{\beta} s_{\beta}x\lambda$ in $\poset{\lambda}$ if and only if $\phi(x)\lambda^\vee \edge{\beta^\vee} s_{\beta^\vee}\phi(x)\lambda^\vee$ in $\poset{\lambda^\vee}$.
        \item Let $\mathfrak{g}_{\af}$ be of type $\mixed$.
              Let $\lambda = \sum_{i \in I}m_i\varpi_i \in  P^0_+$ and $\tilde{\lambda} = \sum_{i \in I}m_i\tilde{\varpi}_i \in \tilde{P}^{0}_{+}$.
              Let $x \in W_\af$ and $\beta \in \positiverealroot$.
              $x\lambda \edge{\beta} s_{\beta}x\lambda$ in $\poset{\lambda}$ if and only if $\psi(x)\tilde{\lambda} \edge{\psi(\beta)} s_{\psi(\beta)}\psi(x)\tilde{\lambda}$ in $\poset{\tilde{\lambda}}$.
    \end{enumerate}
\end{lemma}
\begin{proof}
    We only prove (1); the proof of (2) is similar.
    It is enough to prove the following claim.
    \begin{claim*}
        Let $\lambda = \sum_{i \in I}m_i\varpi_i \in  P^0_+$ and $\lambda^\vee = \sum_{i \in I}m_i\varpi_i^\vee \in P^{0,\vee}_{+}$.
        For $x,y \in W_\af$, $x\lambda \le y\lambda$ in $W_\af$ if and only if $\phi(x)\lambda^\vee \le \phi(y)\lambda^\vee$ in $W_\af^\vee$.
    \end{claim*}
    We prove the ``only if'' part.
    Assume that $x\lambda \le y\lambda$ in $W_\af$.
    By definition, there exists a sequence $\beta_1, \dots , \beta_k$ of positive real roots such that
    \begin{equation*}
        \pairing{\beta_m^\vee}{s_{\beta_{m-1}}\cdots s_{\beta_1}x\lambda} > 0
    \end{equation*}
    for $1 \le m \le k$.
    Hence, we have $\mathrm{cl}(x^{-1}s_{\beta_1} \cdots s_{\beta_{m-1}}\beta_m) \in \Delta^+ \setminus \Delta_{J_\lambda}^+$.
    Since $J_\lambda = J_{\lambda^\vee}$, we see that
    \begin{equation*}
        \pairing{\phi(s_{\beta_{m-1}})\cdots\phi(s_{\beta_1})\phi(x)\lambda^\vee}{\beta_m} > 0
    \end{equation*}
    for $1 \le m \le k$.
    Therefore, we have $\phi(x)\lambda^\vee \le \phi(y)\lambda^\vee$.
    We can prove the ``if'' part in the same way.
\end{proof}

\begin{lemma}\label{lemma:sib}
    \begin{enumerate}
        \item Let $\mathfrak{g}_{\af}$ be of dual untwisted type.
              Let $J \subset I$.
              Let $x \in \pet{J}$ and $\beta \in \positiverealroot$.
              Then, $x \edge{\beta} s_{\beta}x$ in $\sib^J$ if and only if $\phi(x) \edge{\beta^\vee} s_{\beta^\vee}\phi(x)$ in $\sib^{J,\vee}$.
              Here, $\sib^{J,\vee}$ is the parabolic semi-infinite Bruhat graph for $\mathfrak{g}_\af^\vee$.
        \item Let $\mathfrak{g}_{\af}$ be of type $\mixed$.
              Let $J \subset I$.
              Let $x \in \pet{J}$ and $\beta \in \positiverealroot$.
              Then, $x \edge{\beta} s_{\beta}x$ in $\sib^J$ if and only if $\psi(x) \edge{\psi(\beta)} s_{\psi(\beta)}\psi(x)$ in $\frac{\infty}{2}\widetilde{\mathrm{BG}}^{J}$.
              Here, $\frac{\infty}{2}\widetilde{\mathrm{BG}}^{J}$ is the parabolic semi-infinite Bruhat graph for type $\Dtwisted$.
    \end{enumerate}
\end{lemma}
\begin{proof}
    We only prove (1); the proof of (2) is similar.
    We prove the ``only if'' part of this lemma.
    Assume that $x \edge{\beta} s_\beta x$ in $\sib^J$.
    We know that $\sil(x) = \sil(\phi(x))$ for all $x \in W_\af$ from Remark~\ref{remark:semi-length_dual}.
    Thus, we have $\sil(s_{\beta^\vee}\phi(x)) = \sil(\phi(s_\beta x)) = \sil(s_\beta x) = \sil(x) + 1 = \sil(\phi(x)) + 1 $.
    By Lemma~\ref{lemma:coset} and $s_\beta x \in \pet{J}$, we have $\phi(s_\beta x) = s_{\beta^\vee}\phi(x) \in \pet{J}$.
    Hence, we have $\phi(x) \edge{\beta^\vee} s_{\beta^\vee}\phi(x)$ in $\sib^{J,\vee}$.
    We can prove the ``if'' part in the same way.
\end{proof}

From these lemmas, we can derive the comparison lemmas and the result \cite[Lemma 4.1.6]{INS} in the twisted cases.
Therefore, we obtain Lemma \ref{lemma:stability} and Theorem \ref{thm:A}.

\section{Proof of the isomorphism theorem}
\label{section:connected_component}

In this section, we give a proof of Sublemma \ref{slem:condition} for both the cases of dual untwisted types and type $\mixed$, separately.

\subsection{Some technical lemmas}

\begin{lemma}\label{lemma:4.2.2}
    Let $J \subset I$.
    \begin{enumerate}
        \item Let $\mathfrak{g}_\af$ be of dual untwisted type.
              Let $x \in \pet{J}$ and $\beta \in \positiverealroot$ be such that $x \edge{\beta} s_{\beta}x$ in $\sib^J$.
              Then, $\beta$ is either of the following forms: $\beta = \alpha$ with $\alpha \in \Delta^+$, or $\beta = \alpha + c_\alpha\delta$ with $\alpha \in -\Delta^+$.
              Moreover, if $x = wz_{\xi}t_{\xi}$ with $w \in W^J$ and $\xi \in \adjusted{J}$, then $w^{-1}\alpha \in \Delta^+ \setminus \Delta_J^+$ in both cases above.
        \item Let $\mathfrak{g}_\af$ be of type $\mixed$.
              Let $x \in \pet{J}$ and $\beta \in \positiverealroot$ be such that $x \edge{\beta} s_{\beta}x$ in $\sib^J$.
              Then, $\beta$ is either of the following forms: $\beta = \alpha$ with $\alpha \in \Delta^+$, or $\beta = \alpha + \delta$ with a long root $\alpha \in -\Delta^+$, or $\beta = 2\alpha + \delta$ with a short root $\alpha \in -\Delta^+$.
              Moreover, if $x = wz_{\xi}t_{\xi}$ with $w \in W^J$ and $\xi \in \adjusted{J}$, then $w^{-1}\alpha \in \Delta^+ \setminus \Delta_J^+$ in all cases above.
    \end{enumerate}
\end{lemma}
\begin{proof}
    It follows from Lemma \ref{lemma:sib} and \cite[Lemma 4.2.2]{INS}.
\end{proof}

\begin{lemma}\label{lemma:6.1.3}
    Let $J \subset I$.
    \begin{enumerate}
        \item Let $\mathfrak{g}_{\af}$ be of dual untwisted type.
              For all $i\in  J^c:= I\setminus J$ and $\xi \in \adjusted{J}$, there exists a positive real root $\beta \in \positiverealroot$ of the form:
              \begin{equation*}
                  \beta = -\gamma + c_{\gamma}\delta
              \end{equation*}
              with $\gamma \in \Delta^+\setminus\Delta^+_J$ such that $[\gamma]_{J^c}=\alpha_i$, $s_\beta z_\xi t_\xi \in \pet{J}$ and $s_\beta z_\xi t_\xi \edge{\beta} z_\xi t_\xi$ in $\sib^J$.
        \item Let $\mathfrak{g}_{\af}$ be of type $\mixed$.
              For all $i\in  J^c:= I\setminus J$ and $\xi \in \adjusted{J}$, there exists a positive real root $\beta \in \positiverealroot$ of one of the following forms:
              \begin{equation*}
                  \beta = -\gamma + \delta \text{ or } -2\gamma + \delta
              \end{equation*}
              with $\gamma \in \Delta^+\setminus\Delta^+_J$ such that $[\gamma]_{J^c}=\alpha_i$, $s_\beta z_\xi t_\xi \in \pet{J}$ and $s_\beta z_\xi t_\xi \edge{\beta} z_\xi t_\xi$ in $\sib^J$.
    \end{enumerate}
\end{lemma}
\begin{proof}
    Suppose that $\mathfrak{g}_\af$ is of dual untwisted type.
    Through the identification mentioned in \S\ref{subsubsection:dual_untwisted}, we can regard $\xi \in \adjusted{J}$ as a $J$-adjusted element in the coroot lattice of $\mathfrak{g}_\af^{\vee}$.
    By applying \cite[Lemma~6.1.3]{INS} to $i \in J^c$ and $\xi \in \adjusted{J}$, we see that there exist $\beta^\vee = -\gamma^\vee + K \in \Delta_{\af}^{\vee,+}$ with $\gamma \in \Delta^+ \setminus \Delta_J^+$ such that $[\gamma]_{J^c} = \alpha_i$, $s_{\beta^\vee}z_\xi t_\xi \in (W^J)_\af^\vee$, and $s_{\beta^\vee}z_\xi t_\xi \edge{\beta^\vee} z_\xi t_\xi$ in $\sib^{J,\vee}$.
    Therefore, we have $\beta = -\gamma + c_\gamma\delta \in \positiverealroot$, $s_\beta z_\xi t_\xi = \phi^{-1}(s_{\beta^\vee}z_\xi t_\xi) \in \pet{J}$ by Lemma~\ref{lemma:coset}, and $s_\beta z_\xi t_\xi \edge{\beta} z_\xi t_\xi$ in $\sib^J$ by Lemma~\ref{lemma:sib}.

    The proof for the case of type $\mixed$ is reduced to that for the case of type $\Dtwisted$ by using the map $\psi$, Lemma~\ref{lemma:coset}, and Lemma~\ref{lemma:sib}.
\end{proof}

\begin{lemma}\label{lemma:4.1.7}
    Let $J \subset I$, $\xi \in \adjusted{J}$, and $\beta \in \positiverealroot$.
    If $z_\xi t_\xi \edge{\beta} s_\beta z_\xi t_\xi$ in $\sib^J$, then $\beta = \alpha_i$ for some $i \in J^c = I \setminus J$.
\end{lemma}
\begin{proof}
    Suppose that $\mathfrak{g}_\af$ is of dual untwisted type.
    By Lemma~\ref{lemma:sib}, we have $z_{\xi}t_{\xi} \edge{\beta^\vee} s_{\beta^\vee}z_{\xi}t_{\xi}$ in $\sib^{J,\vee}$.
    By applying \cite[Lemma~4.1.7]{INS}, we see that $\beta^{\vee} = \alpha_i^{\vee}$ for some $i \in J^c = I \setminus J$, thus $\beta = \alpha_i$.
    The proof for the case of type $\mixed$ is reduced to that for the case of $\Dtwisted$ by using Lemma~\ref{lemma:sib}.
\end{proof}

\subsection{Proof of Sublemma \ref{slem:condition} for dual untwisted types.}
\label{subsection:proof_slem_dual}

In this subsection, suppose that $\mathfrak{g}_\af$ is of dual untwisted type.

\begin{proof}[Proof of Sublemma \ref{slem:condition}]
    We only need to prove the following fact.
    For $\zeta, \xi \in \adjusted{J_\lambda}$, and $a \in (0,1]_{\mathbb{Q}}$, there exists a nonzero length $a$-chain of shape $\lambda$ from $z_{\zeta}t_{\zeta}$ to $z_{\xi}t_{\xi}$ in $\sib^{J_\lambda}$ if and only if $J^c(\lambda; a) \ne \emptyset$ and $[\xi -\zeta]_{J_{\lambda}^{c}} \in Q^+_{J^c(\lambda;a)} \setminus \set{0}$.

    We prove the ``if'' part.
    For $\alpha = \sum_{i \in I} k_i\alpha_i \in Q^+ \coloneqq \sum_{i \in I} \mathbb{Z}_{\ge 0}\alpha_{i}$, we set $\mathrm{ht}(\alpha) \coloneqq \sum_{i \in I}k_i \geq 0$.
    The proof is by induction on $\mathrm{ht}([\xi-\zeta]_{J_\lambda^c})$.
    If $\mathrm{ht}([\xi-\zeta]_{J_\lambda^c}) = 1$, we have $[\xi-\zeta]_{J_{\lambda}^c} = \alpha_i$ for some $i \in J^c(\lambda;a)$.
    Applying Lemma~\ref{lemma:6.1.3} for $\xi$ and $i$, we see that there exists $\beta \in \positiverealroot$ of the form $\beta = -\gamma + c_\gamma\delta$ with $\gamma \in \Delta^+ \setminus \Delta_{J_\lambda}^+$ such that $[\gamma]_{J_{\lambda}^c} = \alpha_i$, $s_{\beta}z_{\xi}t_{\xi} \in \pet{J_\lambda}$, and $s_{\beta}z_{\xi}t_{\xi} \edge{\beta} z_{\xi}t_{\xi}$ in $\sib^{J_\lambda}$.
    Since $[\gamma]_{J_\lambda^c} = \alpha_i$, we have
    \begin{equation*}
        \pairing{\gamma^\vee}{\lambda} = \pairing{\frac{(\alpha_{i},\alpha_{i})}{(\gamma,\gamma)}\alpha_{i}^{\vee}}{\lambda} = \pairing{\alpha_{i}^{\vee}}{\lambda}\frac{(\alpha_{i},\alpha_{i})}{(\gamma,\gamma)}.
    \end{equation*}
    We also see that $\displaystyle\frac{(\alpha_{i},\alpha_{i})}{(\gamma,\gamma)} \in \mathbb{Z}$ from \cite[Proposition~5.1 d)]{Kac}.
    Since $i \in J^{c}(\lambda;a)$, $a\pairing{\alpha_{i}^{\vee}}{\lambda} \in \mathbb{Z}$, and hence we obtain
    \begin{align*}
        a\pairing{\beta^\vee}{s_{\beta}z_{\xi}t_{\xi}\lambda} & = a\pairing{-\beta^\vee}{\lambda} \quad \text{(since $s_{\beta}\beta = -\beta$ and $z_{\xi}t_{\xi}\lambda = \lambda - (\xi,\lambda)\delta$)} \\
                                                              & = a\pairing{\gamma^\vee}{\lambda} \quad \text{(since $\beta^{\vee} = -\gamma^{\vee} + K$)}                                                   \\
                                                              & = a\pairing{\alpha_i^\vee}{\lambda}\frac{\bilinearform{\alpha_i}{\alpha_i}}{\bilinearform{\gamma}{\gamma}} \in \mathbb{Z}.
    \end{align*}
    Hence, the edge $s_{\beta}z_{\xi}t_{\xi} \edge{\beta} z_{\xi}t_{\xi}$ is an $a$-chain of shape $\lambda$.
    Since $s_\beta = s_\gamma t_{-\gamma}$, we have $s_\beta z_\xi t_\xi = s_\gamma z_\xi t_{\xi-z_\xi^{-1}\gamma} \in \pet{J_\lambda}$ and $\xi-z_\xi^{-1}\gamma \in \adjusted{J_\lambda}$.
    Since $[\xi-z_{\xi}^{-1}\gamma]_{J_\lambda^c} = [\xi]_{J_\lambda^c} - \alpha_i = [\zeta]_{J_\lambda^c}$, we deduce that $\xi-z_{\xi}^{-1}\gamma = \zeta$ by Corollary \ref{corollary:adjusted}.
    Thus, we have $s_\beta z_\xi t_\xi = \lfloor s_\gamma \rfloor^{J_{\lambda}}z_\zeta t_\zeta$ by Lemma \ref{lemma:adjusted}(3).
    Because $[\gamma]_{J_\lambda^c} = \alpha_i$, $\gamma \in \Delta_{I(\lambda;a)}$ if we take reduced expression $\lfloor s_\gamma \rfloor^{J_{\lambda}} = s_{i_1} \cdots s_{i_k}$ of $\lfloor s_\gamma \rfloor^{J_{\lambda}}$, then $i_m \in I(\lambda;a)$ for all $1 \le m \le k$.
    We prove that $s_{i_{m}} \cdots s_{i_{k}} \in W^{J_{\lambda}}$ for $1 \le m \le k$ by reverse induction on $m$.
    Since $\lfloor s_\gamma \rfloor^{J_{\lambda}} \in W^{J_{\lambda}}$, we see that $i_k \notin J_{\lambda}$ and hence $s_{i_{k}} \in W^{J_{\lambda}}$.
    Assume that $s_{i_{m+1}} \cdots s_{i_k} \in W^{J_{\lambda}}$.
    Then, $s_{i_{m+1}} \cdots s_{i_k}\gamma \in \Delta^{+}$ for $\gamma \in \Delta_{J_{\lambda}}^{+}$.
    Since $i_k \notin J_{\lambda}$, we have $s_{i_{k}} \cdots s_{i_{m+1}}\alpha_{i_{m}} \notin \Delta_{J_{\lambda}}$.
    Thus, we obtain $\gamma \neq s_{i_{k}} \cdots s_{i_{m+1}}\alpha_{i_{m}}$ and hence $s_{i_{m+1}} \cdots s_{i_k}\gamma \neq \alpha_{i_{m}}$.
    Then, we have $s_{i_{m}}s_{i_{m+1}} \cdots s_{i_{k}}\gamma \in \Delta^{+}$ and hence $s_{i_{m}}s_{i_{m+1}} \cdots s_{i_{k}} \in W^{J_{\lambda}}$.
    This implies that $s_{i_{m}} \cdots s_{i_{k}}z_{\zeta}t_{\zeta} \in \pet{J_{\lambda}}$ by Lemma \ref{lemma:adjusted}(3).
    Because $\lfloor s_\gamma \rfloor^{J_{\lambda}} = s_{i_1} \cdots s_{i_k}$ is a reduced expression of $\lfloor s_\gamma \rfloor^{J_{\lambda}}$, we have $\sil(s_{i_{m}}s_{i_{m+1}} \cdots s_{i_{k}}z_{\zeta}t_{\zeta}) = \sil(s_{i_{m+1}} \cdots s_{i_{k}}z_{\zeta}t_{\zeta}) + 1$ for $1 \le m \le k$.
    Therefore, we obtain that
    \begin{equation*}
        z_{\zeta}t_{\zeta} \edge{\alpha_{i_k}} s_{i_k}z_{\zeta}t_{\zeta} \edge{\alpha_{i_{k-1}}} \cdots \edge{\alpha_{i_1}} s_{i_1} \cdots s_{i_k}z_{\zeta}t_{\zeta} = \lfloor s_\gamma \rfloor^{J_{\lambda}}z_\zeta t_\zeta = s_{\beta}z_{\xi}t_{\xi}
    \end{equation*}
    in $\sib^{J_\lambda}$.
    Therefore, we obtain an $a$-chain of shape $\lambda$ from $z_\zeta t_\zeta$ to $z_\xi t_\xi$ in $\sib^{J_\lambda}$.
    Assume that $\mathrm{ht}([\xi-\zeta]_{J_\lambda^c}) > 1$.
    We take $i \in J^c(\lambda;a)$ such that $[\xi - \zeta]_{J_\lambda^c} -\alpha_i \in Q_{J^c(\lambda;a)}^+$.
    Set $\xi' \coloneqq \xi-\alpha_i + \phi_{J_\lambda}(\xi-\alpha_i) \in \adjusted{J_\lambda}$.
    Then there exists an $a$-chain of shape $\lambda$ from $z_\zeta t_\zeta$ to $z_{\xi'}t_{\xi'}$ and an $a$-chain of shape $\lambda$ from $z_{\xi'}t_{\xi'}$ to $z_\xi t_\xi$ by our induction hypothesis.
    Therefore, we obtain an $a$-chain of shape $\lambda$ from $z_{\zeta}t_{\zeta}$ to $z_{\xi}t_\xi$ in $\sib^{J_\lambda}$.
    The ``only if'' part can be proved in the same way as \cite[Proposition~6.2.2]{INS}, replacing \cite[Lemma~4.1.7]{INS} with Lemma~\ref{lemma:4.1.7} and \cite[Corollary~4.2.2]{INS} with Lemma~\ref{lemma:4.2.2}.
\end{proof}

\subsection{Proof of Sublemma \ref{slem:condition} for type $\mixed$.}

In this subsection, suppose that $\mathfrak{g}_\af$ is of type $\mixed$.
For $\lambda = \sum_{i \in I} m_i\varpi_i \in P^0_+$ and $a \in (0,1]_{\mathbb{Q}}$, we set
\begin{gather*}
    \overline{J^c(\lambda;a)} \coloneqq \{ i \in J_\lambda^c \mid a\pairing{\alpha_i^\vee}{\lambda} \in \mathbb{Z}\},\\
    \overline{I(\lambda; a)} \coloneqq \set{ i \in I \mid a\pairing{\alpha_i^\vee}{\lambda} \in \mathbb{Z} } = J_{\lambda} \sqcup J^c(\lambda; a).
\end{gather*}

Note that $J^c(\lambda;a) \subset \overline{J^c(\lambda;a)}$ and $I(\lambda;a) \subset \overline{I(\lambda;a)}$ since $\pairing{\alpha_{\ell}^\vee}{\lambda} = 2m_\ell$.

For $\alpha = \sum_{i \in I} k_i\alpha_i \in Q^{+} \coloneqq \sum_{i \in I} \mathbb{Z}_{\ge 0}\alpha_{i}$, we set
\begin{gather*}
    \mathrm{supp}(\alpha) \coloneqq \set{i \in I \mid k_i \neq 0}.
\end{gather*}

\begin{proof}[Proof of Sublemma \ref{slem:condition}]
    We only need to prove the following fact.
    For $\zeta, \xi \in \adjusted{J_\lambda}$, and $a \in (0,1]_{\mathbb{Q}}$, there exists a nonzero length $a$-chain of shape $\lambda$ from $z_{\zeta}t_{\zeta}$ to $z_{\xi}t_{\xi}$ in $\sib^{J_\lambda}$ if and only if $J^c(\lambda; a) \neq \emptyset$ and $[\xi -\zeta]_{J_{\lambda}^c} \in Q^+_{J^c(\lambda;a)} \setminus \set{0}$.

    First, we prove the ``if'' part.
    The proof is by induction on $\mathrm{ht}([\xi-\zeta]_{J_\lambda^c})$.
    If $\mathrm{ht}([\xi-\zeta]_{J_\lambda^c}) = 1$, we have $[\xi-\zeta]_{J_{\lambda}^c} = \alpha_i$ for some $i \in J^c(\lambda;a)$.
    Applying Lemma~\ref{lemma:6.1.3} for $\xi$ and $i$, we see that there exists $\beta \in \positiverealroot$ of one of the following forms $\beta = -\gamma + \delta$ or $-2\gamma + \delta$, with $\gamma \in \Delta^+ \setminus \Delta_{J_\lambda}^+$ such that $[\gamma]_{J_{\lambda}^c} = \alpha_i$, $s_{\beta}z_{\xi}t_{\xi} \in \pet{J_\lambda}$, and $s_{\beta}z_{\xi}t_{\xi} \edge{\beta} z_{\xi}t_{\xi}$ in $\sib^{J_\lambda}$.
    Then, the edge $s_{\beta}z_{\xi}t_{\xi} \edge{\beta} z_{\xi}t_{\xi}$ is an $a$-chain of shape $\lambda$.
    Indeed, if $i \neq \ell$ and $\beta$ is of the form $\beta = -\gamma+\delta$, noting that $[\gamma]_{J_{\lambda}^c} = \alpha_i$, we have
    \begin{equation*}
        a\pairing{\beta^\vee}{s_\beta z_\xi t_\xi} = a\pairing{\gamma^\vee}{\lambda} = a\pairing{\alpha_i^\vee}{\lambda}\frac{(\alpha_i,\alpha_i)}{(\gamma,\gamma)} = am_i \in \mathbb{Z}
    \end{equation*}
    since $\alpha_i$ and $\gamma$ are long roots in $\Delta$.
    If $i \neq \ell$ and $\beta$ is of the form $\beta = -2\gamma +\delta$, we have
    \begin{equation*}
        a\pairing{\beta^\vee}{s_\beta z_\xi t_\xi} = a\pairing{\frac{1}{2}\gamma^\vee}{\lambda} = \frac{1}{2}a\pairing{\alpha_i^\vee}{\lambda}\frac{(\alpha_i,\alpha_i)}{(\gamma,\gamma)} = am_i \in \mathbb{Z}
    \end{equation*}
    since $\alpha_i$ is a long root and $\gamma$ is a short root in $\Delta$.
    If $i = \ell$, since $[\gamma]_{J_{\lambda}^c} = \alpha_{\ell}$ and $\Delta$ is of type $B_{\ell}$, we deduce that $\gamma$ is a short root in $\Delta$.
    Thus, $\beta = -2\gamma + \delta$.
    We have
    \begin{equation*}
        a\pairing{\beta^\vee}{s_\beta z_\xi t_\xi} = a\pairing{\frac{1}{2}\gamma^\vee}{\lambda} = \frac{1}{2}a\pairing{\alpha_\ell^\vee}{\lambda}\frac{(\alpha_\ell,\alpha_\ell)}{(\gamma,\gamma)} = \frac{1}{2}a \cdot 2m_\ell=  am_\ell \in \mathbb{Z}.
    \end{equation*}
    Hence, we see that the edge $s_{\beta}z_{\xi}t_{\xi} \edge{\beta} z_{\xi}t_{\xi}$ is an $a$-chain of shape $\lambda$.
    The rest of the proof of the ``if'' part is the same as that in the case of dual untwisted types.

    Next, we prove the ``only if'' part.
    Let
    \begin{gather*}
        z_\zeta t_\zeta = y_0 \edge{\beta_1} y_1 \edge{\beta_2} \cdots \edge{\beta_k} y_k = z_\xi t_\xi
    \end{gather*}
    be a nonzero length $a$-chain of shape $\lambda$ in $\sib^{J_\lambda}$.
    Since $\zeta \neq \xi$, we deduce that $[\xi - \zeta]_{J_\lambda^c} \neq 0$ by Corollary \ref{corollary:adjusted}.
    Write $y_m$ for $0 \le m \le k$ as
    \begin{equation*}
        y_m = w_{m}z_{\xi_m}t_{\xi_m}, \quad \text{with $w_m \in W^{J_\lambda}$ and $\xi \in \adjusted{J_\lambda}$,}
    \end{equation*}
    by Lemma \ref{lemma:adjusted}, and write $\beta_m$ for $1 \le m \le k$ as one of the following forms
    \begin{equation*}
        \beta_m = w_{m-1}\gamma_m \text{ or } w_{m-1}\gamma_m + \delta \text{ or } 2w_{m-1}\gamma_m + \delta,
    \end{equation*}
    with $\gamma_m \in \Delta^+\setminus \Delta_{J_\lambda}^+$, by Lemma \ref{lemma:4.2.2}(2).
    We set
    \begin{gather*}
        E_B \coloneqq \set{1 \le m \le k \mid \beta_m = w_{m-1}\gamma_m},\\
        E_Q^1 \coloneqq \set{1 \le m \le k \mid \beta_m = w_{m-1}\gamma_m + \delta},\\
        E_Q^2 \coloneqq \set{1 \le m \le k \mid \beta_m = 2w_{m-1}\gamma_m + \delta},\\
        E_Q \coloneqq E_Q^1 \cup E_Q^2.
    \end{gather*}
    Then, we have
    \begin{equation*}
        \xi -\zeta = \sum_{1 \le m \le k} (\xi_m -\xi_{m-1}) = \sum_{m \in E_Q}z_{\xi_{m-1}}^{-1}\gamma_m
    \end{equation*}
    by direct calculation, and hence
    \begin{equation*}
        [\xi - \zeta]_{J_\lambda^c} = \sum_{m \in E_Q}[\gamma_m]_{J_\lambda^c} \in Q_{J_\lambda^c}^+.
    \end{equation*}
    We can prove that $\gamma_m \in \Delta_{\overline{I(\lambda;a)}}$ for $1 \le m \le k$ by the same argument as in the proof of \cite[Proposition~6.2.2]{INS}.
    Thus, we have
    \begin{equation*}
        [\xi - \zeta]_{J_\lambda^c} = \sum_{m \in E_Q}[\gamma_m]_{J_\lambda^c} \in Q_{\overline{J^c(\lambda;a)}}^+.
    \end{equation*}
    When $m_\ell = 0$ or $a \notin \left\{\dfrac{q}{2m_\ell} \, \middle\vert \, 1 \le q \le 2m_\ell, \text{$q$ is odd integer}\right\}$, we see that $\overline{J^c(\lambda;a)} = J^c(\lambda;a)$ and $\overline{I(\lambda;a)} = I(\lambda;a)$, and hence the assertion holds.

    Let us show that $\ell \notin \mathrm{supp}([\xi -\zeta]_{J_\lambda^c})$ when $m_\ell \neq 0$ and $a = \dfrac{q}{2m_\ell}$ for an odd integer $1 \le q \le 2m_\ell$.
    Let us take $m \in E_Q$.
    Note that $\Delta$ is of type $B_\ell$.
    If $m \in E_Q^2$, then $\gamma_m$ is a short root in $\Delta$, $\gamma_m^\vee$ is of the form $\gamma_m^\vee = 2\alpha_i^\vee + \dots + 2\alpha_{\ell-1}^\vee + \alpha_\ell^\vee$ with $i \in I(\lambda;a)$, and $\beta_m^\vee = \dfrac{1}{2}w_{m-1}\gamma_m^\vee + \dfrac{1}{2}K$.
    We see that
    \begin{equation*}
        a\pairing{\beta_m^\vee}{y_{m-1}\lambda} = \frac{q}{2m_\ell}\pairing{\frac{1}{2}\gamma_m^\vee}{\lambda} = \frac{q}{2m_\ell}(m_i + \dots +m_{\ell-1} + m_\ell) \notin \mathbb{Z}
    \end{equation*}
    since $\mathrm{supp}(\gamma_m) \subset \overline{I(\lambda;a)}$.
    This contradicts the assumption that the directed path is an $a$-chain of shape $\lambda$.
    Thus, $m \in E_Q^1$ and $\gamma_m$ is a long root in $\Delta$.
    Since $\ell^{\frac{\infty}{2}}(y_m) = \ell^{\frac{\infty}{2}}(s_{\beta_m}y_{m-1}) = \ell^{\frac{\infty}{2}}(y_{m-1}) + 1$, we deduce that $\ell(w_{m}z_{\xi_m}) = \ell(w_{m-1}z_{\xi_{m-1}}s_{z_{\xi_m}^{-1}\gamma_m}) = \ell(w_{m-1}z_{\xi_{m-1}})+1-2\pairing{\rho^\vee}{z_{\xi_{m-1}}\gamma_m}$ (see Proposition \ref{proposition:siBG_QBG_mixed}).
    Hence, $z_{\xi_{m-1}}^{-1}\gamma_m \in \Delta^+$ should satisfy the condition
    \begin{equation}
        \label{equation:length_condition}
        \ell(s_{z_{\xi_{m-1}}^{-1}\gamma_m}) = 2\pairing{\rho^\vee}{z_{\xi_{m-1}}\gamma_m}-1.
    \end{equation}
    By \cite[\S1.13]{BMO} and \cite[Lemma 7.2]{BMO}, a long root $\alpha \in \Delta^+$ satisfying the condition \eqref{equation:length_condition} is of the form $\alpha = \alpha_i+ \dots + \alpha_{j-1}$ for some $1 \le i < j \le \ell$, and hence we see that $\ell \notin \mathrm{supp}(z_{\xi_{m-1}}^{-1}\gamma_m)$; note that we regard $\Delta^{\vee} \coloneqq \set{\alpha^{\vee} \mid \alpha \in \Delta}$ as a root system of type $C_{\ell}$ and apply the result in \cite{BMO} to $\Delta^{\vee}$ here.
    Since $z_{\xi_{m-1}}^{-1} \in W_{J_{\lambda}}$, it follows that $\mathrm{supp}(\gamma_m) \subset \mathrm{supp}(z_{\xi_{m-1}}^{-1}\gamma_m) \cup J_{\lambda}$.
    The condition $m_{\ell} \neq 0$ implies $\ell \notin J_{\lambda}$ and hence we have $\ell \notin \mathrm{supp}(\gamma_m)$.
    This completes the proof of the sublemma.
\end{proof}

\begin{remark}
    The equality \eqref{equation:length_condition} is the condition for quantum edges of the quantum Bruhat graph.
    See Appendix \ref{appendix:qbg_qls} for details.
\end{remark}

\appendix
\section{The quantum Bruhat graph and the quantum LS paths.}
\label{appendix:qbg_qls}

In this appendix, we describe the quantum Bruhat graph and quantum LS paths, and discuss their relationship with the semi-infinite Bruhat graph and semi-infinite LS paths.

Let $\mathfrak{g}_{\af} = \mathfrak{g}(X_{\af})$ be a twisted affine Lie algebra.
Recall that $X$ is the principal submatrix of $X_{\af}$ indexed by $I$ and $\mathfrak{g} = \mathfrak{g}(X) \subset \mathfrak{g}_{\af}$, which is a finite dimensional simple Lie subalgebra of $\mathfrak{g}_{\af}$.
Let $\overset{\circ}{\varpi}_{i}$ denote the $i$-th fundamental weight of $\mathfrak{g}$ for $i \in I$.
We set
\begin{equation*}
    P^{+} \coloneqq \sum_{i \in I} \mathbb{Z}_{\ge 0} \overset{\circ}{\varpi}_{i} \subset \mathfrak{h}.
\end{equation*}
We see that $\mathrm{cl}(P_{+}^{0}) \subset P^{+}$.
For $\lambda \in P^{+}$, we set $J_{\lambda} \coloneqq \set{i \in I \mid \pairing{\alpha_{i}^{\vee}}{\lambda} = 0}$.
Note that $J_{\lambda} = J_{\mathrm{cl}(\lambda)}$ for $\lambda \in P^{0}_{+}$.

For $J \subset I$, we set $\rho_J^{\vee} \coloneqq \frac{1}{2}\sum_{\alpha \in \Delta_J^+} \alpha^{\vee}$.
We remark that $\pairing{\rho^{\vee}-\rho_J^{\vee}}{\xi} = 0$ for all $\xi \in Q_J$.
We deduce the following lemmas from \cite[Lemma 3.10]{LNSSS1} and \cite[Lemma 2.3.7]{INS}.

\begin{lemma}
    \label{lemma:3.10}
    Let $J \subset I$ and $\xi \in \adjusted{J}$.
    Then, $\ell(z_{\xi}) = -2\pairing{\rho^{\vee}_{J}}{\xi}$.
\end{lemma}

\begin{lemma}
    \label{lemma:2.3.7}
    Let $J \subset I$, $x \in W_{\af}$, and $\xi \in \adjusted{J}$.
    Then, $xz_{\xi}t_{\xi} \in \pet{J}$ if and only if $x \in \pet{J}$.
\end{lemma}

\subsection{QBG and QLS paths for dual untwisted types.}
In this subsection, suppose that $\mathfrak{g}_{\af}$ is of dual untwisted type.
In the case of dual untwisted types, we have $\mathrm{cl}(\varpi_{i}) = \overset{\circ}{\varpi}_{i}$ for $i \in I$ and $\mathrm{cl}(P_{+}^{0}) = P^{+}$.

\begin{definition}[the quantum Bruhat graph]\
    \begin{enumerate}
        \item Let $J \subset I$.
              The (parabolic) quantum Bruhat graph $\qbg^J$ is the $(\Delta^+ \setminus \Delta_J^+)$-labeled directed graph with the vertex set $W^J$. For $w,v  \in W^J$ and $\alpha \in \Delta^+ \setminus\Delta_J^+$,
              \begin{equation*}
                  w \edge{\gamma} v \quad \text{in $\qbg^J$}
              \end{equation*}
              if $v = \lfloor ws_{\gamma} \rfloor^{J}$ and one of the following two conditions holds:
              \begin{enumerate}
                  \item $\ell(v) = \ell(w) + 1$,
                  \item $\ell(v) = \ell(w) - 2\pairing{\rho^{\vee}-\rho_J^{\vee}}{\gamma} + 1$.
              \end{enumerate}
              We call an edge satisfying (a) (resp., (b)) a Bruhat edge (resp., quantum edge).
        \item Let $\lambda \in P^+$ and $a \in (0,1]_{\mathbb{Q}}$. A directed path
              \begin{equation*}
                  w = v_0 \edge{\gamma_1} v_1 \edge{\gamma_2} \cdots \edge{\gamma_k} v_k = v
              \end{equation*}
              from $w$ to $v$ in $\qbg^{J_{\lambda}}$ is an $a$-chain of shape $\lambda$ from $w$ to $v$ if $a\pairing{\gamma_m^{\vee}}{\lambda} \in \mathbb{Z}$ for $1 \le m \le k$.
    \end{enumerate}
\end{definition}

We obtain the following lemma from \cite[Remark 3.1.2]{LNSSS4}, \cite[\S4.3]{LNSSS1} and \cite[\S10]{LS}.

\begin{lemma}
    \label{lemma:QBG_edge_dual}
    Let $J \subset I$.
    \begin{enumerate}
        \item If $w \edge{\gamma} \lfloor ws_{\gamma} \rfloor^{J}$ is a Bruhat edge in $\qbg^{J}$, then $ws_{\gamma} = \lfloor ws_{\gamma} \rfloor^{J} \in W^{J}$.
        \item If $w \edge{\gamma} \lfloor ws_{\gamma} \rfloor^{J}$ is a quantum edge in $\qbg^{J}$, then $\ell(ws_{\gamma}) = \ell(w) - 2\pairing{\rho^{\vee}}{\gamma} + 1$ and $ws_{\gamma}t_{\gamma} \in \pet{J}$.
    \end{enumerate}
\end{lemma}

\begin{definition}[quantum LS paths]
    Let $\lambda \in P^+$.
    $\eta = (w_1,\dots ,w_s;0=a_0<a_1<\dots<a_s=1)$, where $w_m \in W^{J_{\lambda}}$ for $1 \le m \le s$ and $a_m \in (0,1]_{\mathbb{Q}}$ for $1 \le m \le s-1$, is a quantum Lakshmibai--Seshadri path (QLS path) of shape $\lambda$ if $w_{m} \neq w_{m+1}$ and there exists an $a_m$-chain of shape $\lambda$ from $w_{m+1}$ to $w_m$ in $\qbg^{J_{\lambda}}$ for all $1 \le m \le s-1$. $\qls{\lambda}$ denotes the set of all QLS paths of shape $\lambda$.
\end{definition}

\begin{remark}
    Our definition of the quantum Bruhat graph is different from that in \cite[Definition 6.1]{BFP} and \cite[\S4.3]{LNSSS1}, and is the same as \cite[\S4.1]{OS} and \cite[Definition 3.1.2]{N}.
    Note that the edges of the quantum Bruhat graph are labeled by coroots in \cite{N}, whereas they are labeled by roots here.
    To assert the following proposition, we adopt this definition.
\end{remark}

\begin{proposition}
    \label{proposition:siBG_QBG_dual}
    Let $J \subset I$.
    \begin{enumerate}
        \item Let $x = wz_{\xi}t_{\xi} \in \pet{J}$ with $w \in W^{J}$ and $\xi \in \adjusted{J}$ (see Lemma \ref{lemma:adjusted}(3)), $\beta \in \positiverealroot$, and $x \edge{\beta} s_{\beta}x$ in $\sib^{J}$.
              If $\beta = w\gamma$ ($\beta = w\gamma +c_{\gamma}\delta$) with $\gamma \in \Delta^{+} \setminus \Delta_{J}^{+}$, then we have a Bruhat edge (resp., quantum edge) $w \edge{\gamma} \lfloor ws_{\gamma} \rfloor^{J}$ in $\qbg^{J}$.
              (Note Lemma \ref{lemma:4.2.2}(1).)
        \item Let $w \in W^{J}, \gamma \in \Delta^{+} \setminus \Delta_{J}^{+}$, and $w \edge{\gamma} \lfloor ws_{\gamma} \rfloor^{J}$ in $\qbg^{J}$.
              If this edge is a Bruhat edge (resp., quantum edge), then we have $wz_{\xi}t_{\xi} \edge{\beta} s_{\beta}wz_{\xi}t_{\xi}$ in $\sib^{J}$ for all $\xi \in \adjusted{J}$, where $\beta = w\gamma \in \positiverealroot$ (resp., $\beta = w\gamma + c_{\gamma}\delta \in \positiverealroot$) .
    \end{enumerate}
\end{proposition}
\begin{proof}
    (1) From the assumption, we have
    \begin{equation}
        \label{equation:proof}
        \sil(s_{\beta}x) = \sil(x) + 1.
    \end{equation}
    Let $\beta = w\gamma$ with $\gamma \in \Delta^{+} \setminus \Delta_{J}^{+}$.
    Then, $s_{\beta}x = s_{w\gamma}wz_{\xi}t_{\xi} = ws_{\gamma}z_{\xi}t_{\xi}$.
    Substituting $s_{\beta}x = ws_{\gamma}z_{\xi}t_{\xi}$ into \eqref{equation:proof}, we obtain that
    \begin{equation*}
        \ell(ws_{\gamma}z_{\xi}) = \ell(wz_{\xi}) + 1.
    \end{equation*}
    Since $s_{\beta}x \in \pet{J}$, we deduce that $ws_{\gamma} \in W^{J}$ by Lemma \ref{lemma:adjusted}(3).
    Hence, we have $\ell(ws_{\gamma}z_{\xi}) = \ell(ws_{\gamma}) + \ell(z_{\xi})$ and $\ell(wz_{\xi}) = \ell(w) + \ell(z_{\xi})$.
    Therefore, it follows that $\ell(ws_{\gamma}) = \ell(w) + 1$.

    Let $\beta = w\gamma + c_{\gamma}\delta$ with $\gamma \in \Delta^{+}\setminus\Delta_{J}^{+}$.
    Then, $s_{\beta}x = s_{w\gamma}t_{w\gamma}wz_{\xi}t_{\xi} = ws_{\gamma}z_{\xi}t_{\xi+z_{\xi}^{-1}\gamma}$ (note \eqref{equation:semiproduct_dual}).
    Since $s_{\beta}x \in \pet{J}$, we deduce that $\xi+z_{\xi}^{-1}\gamma \in \adjusted{J}$ and $s_{\beta}x = \lfloor ws_{\gamma} \rfloor^{J}z_{\xi+z_{\xi}^{-1}\gamma}t_{\xi+z_{\xi}^{-1}\gamma}$ by Lemma \ref{lemma:adjusted}(3).
    Substituting this into \eqref{equation:proof}, we obtain that
    \begin{equation*}
        \ell(\lfloor ws_{\gamma} \rfloor^{J}) + \ell(z_{\xi+z_{\xi}^{-1}\gamma}) + 2\pairing{\rho^{\vee}}{\xi+z_{\xi}^{-1}\gamma} = \ell(w) + \ell(z_{\xi}) +2\pairing{\rho^{\vee}}{\xi} + 1.
    \end{equation*}
    By Lemma \ref{lemma:3.10}, we have $\ell(\lfloor ws_{\gamma} \rfloor^{J}) = \ell(w) -2\pairing{\rho^{\vee} - \rho_{J}^{\vee}}{\gamma} + 1$.
    Note that $\pairing{\rho^{\vee}-{\rho}_{J}^{\vee}}{z_{\xi}^{-1}\gamma} = \pairing{\rho^{\vee} - \rho_{J}^{\vee}}{\gamma}$ because $z_{\xi}^{-1} \in W_{J}$ and $\pairing{\rho^{\vee}-\rho_J^{\vee}}{\zeta} = 0$ for all $\zeta \in Q_J$.

    (2) Suppose that $w \edge{\gamma} \lfloor ws_{\gamma} \rfloor^{J}$ is a Bruhat edge in $\qbg^{J}$.
    Since $\ell(ws_{\gamma}) > \ell(w)$, we deduce that $\beta = w\gamma \in \Delta^{+} \subset \positiverealroot$ from \cite[Proposition 4.4.6]{BB}.
    By Lemma \ref{lemma:QBG_edge_dual}(1) and Lemma \ref{lemma:adjusted}(3), we see that $ws_{\gamma} \in W^{J}$ and $s_{\beta}wz_{\xi}t_{\xi} = ws_{\gamma}z_{\xi}t_{\xi} \in \pet{J}$.
    Therefore, we have
    \begin{align*}
        \sil(s_{\beta}wz_{\xi}t_{\xi}) & = \sil(ws_{\gamma}z_{\xi}t_{\xi})                                                                              \\
                                       & = \ell(ws_{\gamma}) + \ell(z_{\xi}) + 2\pairing{\rho^{\vee}}{\xi} \quad \text{(since $ws_{\gamma} \in W^{J}$)} \\
                                       & =\ell(w) + 1 + \ell(z_{\xi}) + 2\pairing{\rho^{\vee}}{\xi}                                                     \\
                                       & \quad \text{(since $w \edge{\gamma} \lfloor ws_{\gamma} \rfloor^{J}$ is a Bruhat edge in $\qbg^{J}$)}          \\
                                       & = \ell(wz_{\xi}) + 2\pairing{\rho^{\vee}}{\xi} + 1 = \sil(wz_{\xi}t_{\xi}) + 1.
    \end{align*}

    Suppose that $w \edge{\gamma} \lfloor ws_{\gamma} \rfloor^{J}$ is a quantum edge in $\qbg^{J}$.
    Set $\beta = w\gamma + c_{\gamma}\delta \in \positiverealroot$.
    By Lemma \ref{lemma:QBG_edge_dual}(2) and Lemma \ref{lemma:2.3.7}, we have $ws_{\gamma}t_{\gamma} \in \pet{J}$ and $s_{\beta}wz_{\xi}t_{\xi} = ws_{\gamma}t_{\gamma}z_{\xi}t_{\xi} \in \pet{J}$ (note \eqref{equation:semiproduct_dual}).
    By Lemma \ref{lemma:adjusted}(3), we deduce that $\xi + z_{\xi}^{-1}\gamma \in \adjusted{J}$, and $s_{\beta}wz_{\xi}t_{\xi} = \lfloor ws_{\gamma} \rfloor^{J}z_{\xi + z_{\xi}^{-1}\gamma}t_{\xi + z_{\xi}^{-1}\gamma}$.
    Therefore, we see that
    \begin{align*}
        \sil(s_{\beta}wz_{\xi}t_{\xi}) & = \sil(\lfloor ws_{\gamma} \rfloor^{J}z_{\xi + z_{\xi}^{-1}\gamma}t_{\xi + z_{\xi}^{-1}\gamma})                                    \\
                                       & = \ell(\lfloor ws_{\gamma} \rfloor^{J}) + \ell(z_{\xi + z_{\xi}^{-1}\gamma}) + 2\pairing{\rho^{\vee}}{\xi + z_{\xi}^{-1}\gamma}    \\
                                       & = \ell(w) -2\pairing{\rho^{\vee}-\rho_{J}^{\vee}}{\gamma} + 1                                                                      \\
                                       & \phantom{{}={}} -2\pairing{\rho_{J}^{\vee}}{\xi +z_{\xi}^{-1}\gamma}                                                               \\
                                       & \phantom{{}={}} +2\pairing{\rho^{\vee}}{\xi} + 2\pairing{\rho^{\vee}}{z_{\xi}^{-1}\gamma}                                          \\
                                       & \quad \text{(since $w \edge{\gamma} \lfloor ws_{\gamma} \rfloor^{J}$ is a quantum edge in $\qbg^{J}$}                              \\
                                       & \quad \quad \text {and Lemma \ref{lemma:3.10})}                                                                                    \\
                                       & = \ell(w) -2\pairing{\rho^{\vee}-\rho_{J}^{\vee}}{\gamma} + 1                                                                      \\
                                       & \phantom{{}={}} + \ell(z_{\xi}) -2\pairing{\rho_{J}^{\vee}}{z_{\xi}^{-1}\gamma}                                                    \\
                                       & \phantom{{}={}} +2\pairing{\rho^{\vee}}{\xi} + 2\pairing{\rho^{\vee}}{z_{\xi}^{-1}\gamma} \quad \text{(by Lemma \ref{lemma:3.10})} \\
                                       & = \ell(w) + 1 + \ell(z_{\xi}) + 2\pairing{\rho^{\vee}}{\xi}                                                                        \\
                                       & \phantom{{}={}} + 2\pairing{\rho^{\vee}-\rho_{J}^{\vee}}{z_{\xi}^{-1}\gamma} - 2\pairing{\rho^{\vee}-\rho_{J}^{\vee}}{\gamma}      \\
                                       & = \sil(wz_{\xi}t_{\xi}) + 1.
    \end{align*}
    The final equality is obtained by using $z_{\xi}^{-1} \in W_{J}$ and $\pairing{\rho^{\vee}-\rho_J^{\vee}}{\zeta} = 0$ for all $\zeta \in Q_J$.
\end{proof}

\begin{corollary}
    \label{corollary:siLS_QLS_dual}
    Let $\lambda \in P^{0}_{+}.$
    We have a map $\sils{\lambda} \rightarrow \qls{\mathrm{cl}(\lambda)}$, given by
    \begin{equation*}
        (x_{1},\dots,x_{s};a_{0},\dots.a_{s}) \mapsto (w_{1},\dots,w_{s};a_{0},\dots,a_{s}),
    \end{equation*}
    where $x_{m} = w_{m}z_{\xi_{m}}t_{\xi_{m}}$ with $w_{m} \in W^{J_{\lambda}}$ and $\xi_{m} \in \adjusted{J_{\lambda}}$ for $1 \le m \le s$.
    Moreover, this map is surjective.
\end{corollary}
\begin{proof}
    Let $x = wz_{\xi}t_{\xi}, y = vz_{\zeta}t_{\zeta} \in \pet{J}$ with $w,v\in W^{J_{\lambda}}$ and $\xi,\zeta\in \adjusted{J_{\lambda}}$, and $a \in (0,1]_{\mathbb{Q}}$.
    Assume that there exists an $a$-chain of shape $\lambda$
    \begin{equation*}
        x = y_0 \edge{\beta_{1}} y_{1} \edge{\beta_2} \cdots \edge{\beta_{k}} y_{k} = y
    \end{equation*}
    in $\sib^{J_{\lambda}}$.
    We write $y_{n} \in \pet{J_{\lambda}}$ ($0 \le n \le k$) as $y_{n} = v_{n}z_{\zeta_{n}}t_{\zeta_{n}}$ with $v_{n} \in W^{J_{\lambda}}$ and $\zeta_{n} \in \adjusted{J_{\lambda}}$, and write $\beta_{n}$ ($1 \le n \le k$) as $\beta_{n} = v_{n-1}\gamma_{n} + p_{n}c_{\gamma}\delta$ with $\gamma_{n} \in \Delta^{+} \setminus \Delta_{J_{\lambda}^{+}}$ and $p_{n} \in \set{0,1}$ (see Lemma \ref{lemma:adjusted}(3) and Lemma \ref{lemma:4.2.2}(1)).
    Then, by Proposition \ref{proposition:siBG_QBG_dual}(1), we obtain a directed path
    \begin{equation*}
        w = v_{0} \edge{\gamma_{1}} v_{1} \edge{\gamma_{2}} \cdots \edge{\gamma_{k}} v_{n} = v
    \end{equation*}
    in $\qbg^{J_{\mathrm{cl}(\lambda)}}$.
    By the definition of $a$-chains of shape $\lambda$ in $\sib^{J_{\lambda}}$, we have $a\pairing{\beta_{n}^{\vee}}{y_{n-1}\lambda} \in \mathbb{Z}$ for $1 \le n \le k$.
    Noting $\beta_{n}^{\vee} = v_{n-1}\gamma_{n}^{\vee} + p_{n}K$, we see that
    \begin{align*}
        a\pairing{\beta_{n}^{\vee}}{y_{n-1}\lambda} & = a\pairing{v_{n-1}\gamma_{n}^{\vee} + p_{n}K}{v_{n-1}z_{\zeta_{n-1}}t_{\zeta_{n-1}}\lambda}                 \\
                                                    & = a\pairing{v_{n-1}\gamma_{n}^{\vee}}{v_{n-1}\lambda}                                                        \\
                                                    & = a\pairing{\gamma_{n}^{\vee}}{\lambda} = a\pairing{\gamma_{n}^{\vee}}{\mathrm{cl}(\lambda)} \in \mathbb{Z}.
    \end{align*}
    Hence, the directed path in $\qbg^{J_{\mathrm{cl}(\lambda)}}$ is an $a$-chain of shape $\mathrm{cl}(\lambda)$.
    Therefore, we obtain the map $\sils{\lambda} \rightarrow \qls{\mathrm{cl}(\lambda)}$.
    We prove that this map is surjective.
    Let $\eta = (w_{1},\dots,w_{s};a_{0},\dots,a_{s}) \in \qls{\mathrm{cl}(\lambda)}$.
    We need to find elements $\xi_{1},\dots,\xi_{s} \in \adjusted{J_{\lambda}}$ such that $(w_{1}z_{\xi_{1}}t_{\xi_{1}},\dots,w_{s}z_{\xi_{s}}t_{\xi_{s}};a_{0},\dots,a_{s}) \in \sils{\lambda}$.
    They are constructed by reverse induction on $1 \le m \le s$.
    For $m = s$, set $\xi_{s} = 0$.
    For $m < s$, assume that there are elements $\xi_{m+1},\dots,\xi_{s-1},\xi_{s}(=0) \in \adjusted{J_{\lambda}}$ such that there exists an $a_{q}$-chain of shape $\lambda$ from $w_{q+1}z_{\xi_{q+1}}t_{\xi_{q+1}}$ to $w_{q}z_{\xi_{q}}t_{\xi_{q}}$ in $\sib^{J_{\lambda}}$ for each $m+1 \le q \le s$.
    Since $\eta \in \qls{\mathrm{cl}(\lambda)}$, there exists an $a_{m}$-chain of shape $\mathrm{cl}(\lambda)$
    \begin{equation}
        \label{equation:chain_QBG}
        w_{m+1} = v_{0} \edge{\gamma_{1}} v_{1} \edge{\gamma_{2}} \cdots \edge{\gamma_{k}} v_{k} = w_{m}
    \end{equation}
    in $\qbg^{J_{\mathrm{cl}(\lambda)}}$.
    We define $\beta_{n} \in \Delta_{\af}^{+}$ for $1 \le n \le k$ by
    \begin{equation*}
        \beta_{n} \coloneqq
        \begin{cases}
            v_{n-1}\gamma_{n}                        & \text{if $v_{n-1}\edge{\gamma_{n}}v_{n}$ is a Bruhat edge,}  \\
            v_{n-1}\gamma_{n} + c_{\gamma_{n}}\delta & \text{if $v_{n-1}\edge{\gamma_{n}}v_{n}$ is a quantum edge,}
        \end{cases}
    \end{equation*}
    and define $\zeta_{n} \in \adjusted{J_{\lambda}}$ for $0 \le n \le k$ inductively by
    \begin{equation*}
        \zeta_{n} \coloneqq
        \begin{cases}
            \xi_{m+1}                                    & \text{if $n=0$,}                                                         \\
            \zeta_{n-1}                                  & \text{if $n > 0$ and $v_{n-1}\edge{\gamma_{n}}v_{n}$ is a Bruhat edge,}  \\
            \zeta_{n-1} + z_{\zeta_{n-1}}^{-1}\gamma_{n} & \text{if $n > 0$ and $v_{n-1}\edge{\gamma_{n}}v_{n}$ is a quantum edge.} \\
        \end{cases}
    \end{equation*}
    Then, by Lemma \ref{proposition:siBG_QBG_dual}(2), we obtain a directed path
    \begin{equation*}
        w_{m+1}z_{\xi_{m+1}}t_{\xi_{m+1}} = v_{0}z_{\zeta_{0}}t_{\zeta_{0}} \edge{\beta_{1}} v_{1}z_{\zeta_{1}}t_{\zeta_{1}} \edge{\beta_{2}} \cdots \edge{\beta_{k}} v_{k}z_{\zeta_{k}}t_{\zeta_{k}} = w_{m}z_{\zeta_{k}}t_{\zeta_{k}}
    \end{equation*}
    in $\sib^{J_{\lambda}}$.
    Since \eqref{equation:chain_QBG} is an $a_{m}$-chain of shape $\mathrm{cl}(\lambda)$ in $\qbg^{J_{\lambda}}$, we deduce that this is an $a_{m}$-chain of shape $\lambda$ in $\sib^{J_{\lambda}}$.
    Therefore, $\zeta_{k} \in \adjusted{J_{\lambda}}$ is the desired element for $\xi_{m}$.
    This completes the proof of the corollary.
\end{proof}

\begin{remark}
    The proofs of Proposition \ref{proposition:siBG_QBG_dual} and Corollary \ref{corollary:siLS_QLS_dual} are identical to that given by Nomoto \cite{N} for type $\mixed$ (see the proofs of \cite[Proposition A.3.5]{N} and \cite[Theorem A.2.5]{N}).
    However, we present the proof for general twisted affine types here.
\end{remark}

\subsection{QBG and QLS paths for type $\mixed$.}

In this appendix, suppose that $\mathfrak{g}_{\af}$ is of type $\mixed$.
In this case, we have
\begin{equation*}
    \mathrm{cl}(\varpi_{i}) =
    \begin{cases}
        \overset{\circ}{\varpi}_{i}     & \text{if $i \neq \ell$,} \\
        2\overset{\circ}{\varpi}_{\ell} & \text{if $i = \ell$,}
    \end{cases}
\end{equation*}
for $i \in I$ and $\mathrm{cl}(P^{0}_{+}) \subsetneq P^{+}$.

\begin{definition}
    \label{definition:qbg_mixed}\
    \begin{enumerate}
        \item Let $J \subset I$. The parabolic quantum Bruhat graph $\qbg_{\mixed}^{J}$ is the $(\Delta^+ \setminus \Delta_J^+)$-labeled directed graph with the vertex set $W^J$. For $w,v  \in W^J$ and $\alpha \in \Delta^+ \setminus\Delta_J^+$,
              \begin{equation*}
                  w \edge{\alpha} v \quad \text{in $\qbg_{\mixed}^{J}$}
              \end{equation*}
              if $v = \lfloor ws_\alpha \rfloor^{J}$ and one of the following two conditions holds:
              \begin{enumerate}
                  \item\label{Bruhat} $\ell(v) = \ell(w) + 1$,
                  \item\label{quantum} $\ell(v) = \ell(w) - 2\pairing{\rho^{\vee} - \rho_{J}^{\vee}}{\alpha} + 1$.
              \end{enumerate}
              We call an edge satisfying (a) (resp., (b)) a Bruhat edge (resp., quantum edge).
        \item Let $\lambda \in P^+$ and $a \in (0,1]_{\mathbb{Q}}$. A directed path
              \begin{equation*}
                  w = v_0 \edge{\gamma_1} v_1 \edge{\gamma_2} \cdots \edge{\gamma_k} v_k = v
              \end{equation*}
              from $w$ to $v$ in $\qbg_{\mixed}^J$ is an $a$-chain of shape $\lambda$ from $w$ to $v$ if the edges $v_{m-1} \edge{\gamma_m} v_m$ for $1 \le m \le k$ satisfy
              the following conditions:\\
              If the edge is a Bruhat edge, then $a\pairing{\gamma_{m}^{\vee}}{\lambda} \in \mathbb{Z}$.
              If the edge is a quantum edge, then
              \begin{equation*}
                  \begin{cases}
                      a\pairing{\gamma_{m}^{\vee}}{\lambda} \in \mathbb{Z}  & \text{if $\gamma_{m}$ is a long root in $\Delta$,}  \\
                      a\pairing{\gamma_{m}^{\vee}}{\lambda} \in 2\mathbb{Z} & \text{if $\gamma_{m}$ is a short root in $\Delta$.}
                  \end{cases}
              \end{equation*}
    \end{enumerate}
\end{definition}

\begin{remark}
    \label{remark:QBG_edge_mixed}
    Lemma \ref{lemma:QBG_edge_dual} also holds for $\qbg_{\mixed}^{J}$.
\end{remark}

\begin{definition}[$\mixed$-type QLS paths]
    Let $\lambda \in P^+$.
    $\eta = (w_1,\dots ,w_s;0=a_0<a_1<\dots<a_s=1)$, where $w_m \in W^{J_{\lambda}}$ for $1 \le m \le s$ and $a_m \in (0,1]_{\mathbb{Q}}$ for $1 \le m \le s-1$, is an $\mixed$-type QLS path (or a QLS path of type $\mixed$) of shape $\lambda$ if $w_{m} \neq w_{m+1}$ and there exists an $a_m$-chain of shape $\lambda$ from $w_{m+1}$ to $w_m$ in $\qbg_{\mixed}^{J_{\lambda}}$ for all $1 \le m \le s-1$.
    Let $\qlsmixed{\lambda}$ denote the set of all $\mixed$-type QLS paths of shape $\lambda$.
\end{definition}

\begin{remark}
    In \cite{N}, Nomoto introduced two variants of the graph: $(\mathrm{QBG}_{b\lambda}^{A_{2n}^{(2)}})^{S}$ (\cite[Definition 3.1.4]{N}) and $((\mathrm{QBG}_{b\lambda}^{\dagger})^{A_{2n}^{(2)}})^{S}$ (\cite[Definition 3.2.3]{N}), and defined the $A_{2\ell}^{(2)}$-type QLS paths.
    Our definition of the $a$-chains in $\mathrm{QBG}_{\mixed}^{J}$ differs from the former in that the parity condition for short roots (i.e., belonging to $2\mathbb{Z}$) is imposed on quantum edges, whereas in \cite{N}, it is imposed on Bruhat edges.
    Our definition is also distinct from the latter, which involves half-integer values ($\frac{1}{2}\mathbb{Z}$).
\end{remark}

\begin{proposition}
    \label{proposition:siBG_QBG_mixed}
    Let $J \subset I$.
    \begin{enumerate}
        \item Let $x = wz_{\xi}t_{\xi} \in \pet{J}$ with $w \in W^{J}$ and $\xi \in \adjusted{J}$ (see Lemma \ref{lemma:adjusted}(3)), $\beta \in \positiverealroot$, and $x \edge{\beta} s_{\beta}x$ in $\sib^{J}$.
              If $\beta = w\gamma$ ($\beta = w\gamma + \delta$ or $2w\gamma + \delta$) with $\gamma \in \Delta^{+} \setminus \Delta_{J}^{+}$, then we have a Bruhat edge (resp., quantum edge) $w \edge{\gamma} \lfloor ws_{\gamma} \rfloor^{J}$ in $\qbg_{\mixed}^{J}$ (note Lemma \ref{lemma:4.2.2}(2)).
        \item Let $w \in W^{J}, \gamma \in \Delta^{+} \setminus \Delta_{J}^{+}$, and $w \edge{\gamma} \lfloor ws_{\gamma} \rfloor^{J}$ in $\qbg_{\mixed}^{J}$.
              If this edge is a Bruhat edge, then we have $wz_{\xi}t_{\xi} \edge{\beta} s_{\beta}wz_{\xi}t_{\xi}$ in $\sib^{J}$ for all $\xi \in \adjusted{J}$, where $\beta = w\gamma \in \positiverealroot$.
              If this edge is a quantum edge and $\gamma$ is a long root (short root) in $\Delta$, then we have $wz_{\xi}t_{\xi} \edge{\beta} s_{\beta}wz_{\xi}t_{\xi}$ in $\sib^{J}$ for all $\xi \in \adjusted{J}$, where $\beta = w\gamma + \delta \in \positiverealroot$ (resp., $\beta = 2w\gamma + \delta \in \positiverealroot$).
    \end{enumerate}
\end{proposition}
\begin{proof}
    The proof is the same as that of Proposition \ref{proposition:siBG_QBG_dual}.
    Note \eqref{equation:semiproduct_mixed1}, \eqref{equation:semiproduct_mixed2}, and Remark \ref{remark:QBG_edge_mixed}.
\end{proof}

\begin{corollary}
    Let $\lambda \in P^{0}_{+}.$
    We have a map $\sils{\lambda} \rightarrow \qlsmixed{\mathrm{cl}(\lambda)}$, given by
    \begin{equation*}
        (x_{1},\dots,x_{s};a_{0},\dots.a_{s}) \mapsto (w_{1},\dots,w_{s};a_{0},\dots,a_{s}),
    \end{equation*}
    where $x_{m} = w_{m}z_{\xi_{m}}t_{\xi_{m}}$ with $w_{m} \in W^{J_{\lambda}}$ and $\xi_{m} \in \adjusted{J_{\lambda}}$ for $1 \le m \le s$.
    Moreover, this map is surjective.
\end{corollary}
\begin{proof}
    The proof is the same as that of Corollary \ref{corollary:siLS_QLS_dual}, under the following modification.
    Let $x = wz_{\xi}t_{\xi}$ with $w \in W^{J_{\lambda}}$ and $\xi \in \adjusted{J_{\lambda}}$, $\beta = 2w\gamma + \delta \in \positiverealroot$ with $\gamma \in \Delta^{+} \setminus \Delta_{J_{\lambda}}^{+}$, and $x \edge{\beta} s_{\beta}x$ in $\sib^{J_{\lambda}}$.
    Note that $\gamma$ is a short root in $\Delta$.
    Then, by Proposition \ref{proposition:siBG_QBG_mixed}(1), we obtain a quantum edge $w \edge{\gamma} \lfloor ws_{\gamma}\rfloor^{J_{\mathrm{cl}(\lambda)}}$ in $\qbg_{\mixed}^{J_{\mathrm{cl}(\lambda)}}$.
    If $x \edge{\beta} s_{\beta}x$ is an $a$-chain of shape $\lambda$ for some $a \in (0,1]_{\mathbb{Q}}$, then
    \begin{align*}
        a\pairing{\beta^{\vee}}{x\lambda} & = a\pairing{\frac{1}{2}w\gamma^{\vee}+\frac{1}{2}K}{wz_{\xi}t_{\xi}\lambda}                                                \\
                                          & = \frac{1}{2}a\pairing{\gamma^{\vee}}{\lambda} = \frac{1}{2}a\pairing{\gamma^{\vee}}{\mathrm{cl}(\lambda)} \in \mathbb{Z}.
    \end{align*}
    (Note $\beta^{\vee} = \frac{1}{2}w\gamma^{\vee}+\frac{1}{2}K$.)
    Thus, we have $a \pairing{\gamma^{\vee}}{\mathrm{cl}(\lambda)} \in 2\mathbb{Z}$, and the quantum edge $w \edge{\gamma} \lfloor ws_{\gamma}\rfloor^{J_{\mathrm{cl}(\lambda)}}$ labeled by the short root $\gamma$ is an $a$-chain of shape $\mathrm{cl}(\lambda)$ in $\qbg_{\mixed}^{J_{\mathrm{cl}(\lambda)}}$.
    Conversely, if a quantum edge $w \edge{\gamma} \lfloor ws_{\gamma}\rfloor^{J_{\mathrm{cl}(\lambda)}}$ in $\qbg_{\mixed}^{J_{\mathrm{cl}(\lambda)}}$ with $w \in W^{J_{\mathrm{cl}(\lambda)}}$ and a short root $\gamma \in \Delta \setminus \Delta_{J_{\mathrm{cl}(\lambda)}}^{+}$ is an $a$-chain of shape $\mathrm{cl}(\lambda)$ for some $a \in (0,1]_{\mathbb{Q}}$, then we deduce that $wz_{\xi}t_{\xi} \edge{2w\gamma + \delta} \lfloor ws_{\gamma} \rfloor^{J_{\lambda}}z_{\xi+z_{\xi}^{-1}\gamma}t_{\xi+z_{\xi}^{-1}\gamma}$ is an $a$-chain of shape $\lambda$ in $\sib^{J_{\lambda}}$ for all $\xi \in \adjusted{J_{\lambda}}$ by Proposition \ref{proposition:siBG_QBG_mixed}(2).
\end{proof}

\section*{Acknowledgments.}
The authors are grateful to their supervisor, Takeshi Ikeda, for his helpful guidance throughout this research and for his careful comments on the manuscript.
The authors also wish to thank Satoshi Naito for his valuable comments on the manuscript and stimulating discussions.
Finally, the authors thank Takafumi Kouno for helpful discussions.

\bibliographystyle{alpha}
\bibliography{ref.bib}

\end{document}